\documentclass[10pt,a4paper,reqno, oneside]{amsart}
\usepackage[foot]{amsaddr}
\usepackage{graphicx}

\makeatletter
\renewcommand\subsection{\@startsection{subsection}{1}
	\z@{.5\linespacing\@plus.7\linespacing}{-.5em}
	{\normalfont\scshape}}
\makeatother

\RequirePackage[utf8]{inputenc}
\RequirePackage[T1]{fontenc}
\RequirePackage[english]{babel}
\RequirePackage[toc, page]{appendix}

\RequirePackage[left=3cm,right=5cm,top=4cm,bottom=4cm]{geometry}
\RequirePackage{amsmath, bbm, amsfonts, amssymb, amsthm} % packages for mathematics
\RequirePackage{hyperref, verbatim}
\RequirePackage{enumerate}  % enumerate with custom symbols
\RequirePackage{xcolor, extarrows}
\RequirePackage{graphicx, tikz, tikz-cd}
\RequirePackage[autostyle=true]{csquotes} % for quotation
\RequirePackage{aligned-overset}
\RequirePackage{breakcites}
\RequirePackage{breqn}
\RequirePackage{url}
\RequirePackage{enumitem}
\usepackage{nomencl}
\newcommand{\set}[1]{\left\{ #1 \right\}} % set with bracket notation
\newcommand{\setv}[2]{\left\{ #1 \,\middle|\, #2 \right\}} % set with bracket and vertical line
\newcommand{\setc}[2]{\left\{ #1 : \, #2 \right\}} % set with bracket and colon
\newcommand{\kl}{\left(} % left round bracket
\newcommand{\kr}{\right)} % right round bracket
\newcommand{\abs}[1]{\left| #1 \right|} % absolute value in one command
\newcommand{\N}{\mathbb{N}} % natural numbers
\newcommand{\Z}{\mathbb{Z}} % whole numbers
\newcommand{\R}{\mathbb{R}} % real numbers
\newcommand{\C}{\mathbb{C}} % complex numbers
\newcommand{\symdiff}{\bigtriangleup} % symmetric difference

\newcommand{\Per}{\mathrm{Per}}	% Periodic positions of Toeplitz sequence
\newcommand{\ap}{\mathrm{ap}} % frequency of appearances of blocks
\newcommand{\diag}{\mathrm{diag}} % diagonal matrix

\newcommand{\tm}{\subseteq}
\newcommand{\mt}{\supseteq}

\newcommand{\htop}{h_{\mathrm{top}}}

\newtheorem{thm}{Theorem}[section]
\newtheorem{lem}[thm]{Lemma}
\newtheorem{prop}[thm]{Proposition}

\newtheorem{cor}[thm]{Corollary}

\newtheorem*{thm*}{Theorem}

\newtheoremstyle{named}{}{}{\itshape}{}{\bfseries}{.}{.5em}{\thmnote{#3}}
\theoremstyle{named}
\newtheorem*{namedtheorem}{Theorem}

\theoremstyle{definition}

\newtheorem{defn}[thm]{Definition}
\newtheorem{rem}[thm]{Remark}
\newtheorem{ex}[thm]{Example}

\newcommand{\1}{\mathbbmss{1}} % 1 with double stroke

\newcommand{\ol}[1]{\overline{#1}} % overline

\title{Entropy and strict ergodicity for multidimensional Toeplitz subshifts}
\author{Jamal Drewlo}
\address{Institute of Mathematics, Friedrich Schiller University Jena, 07743 Jena, Germany}
\email{jamal.drewlo@uni-jena.de}

\begin{document}

\begin{abstract}
In this work, we present a self-contained construction that proves the existence of strictly ergodic Toeplitz $\Z^d$-subshifts with arbitrary prescribed entropy.
Moreover, any of these constructed subshifts will have the same maximal equicontinuous factor, which is the universal odometer over $\Z^d$. \\[2pt]
\noindent{\em 2010 Mathematics Subject Classification.} 37B10 (primary), 37B40, 37B05 (secondary).
\end{abstract}

\maketitle

\section{Introduction}
Toeplitz sequences and their associated Toeplitz flows have been introduced to the study of dynamical systems by Jacobs and Keane in \cite{JacobsKeane1969Toeplitz}.
In the subsequent years, their dynamical properties have been thoroughly studied (see \cite{Williams1984ToeplitzFlows}, \cite{DownarowiczLacroix1998Almost1-1} or \cite{Downarowicz2005Survey} for an overview) and they have proved to be very useful for providing examples of dynamical systems with certain prescribed properties (see for example \cite{IwanikLacroix1994NonRegular}, \cite{Iwanik1996ToeplitzWithPP}, \cite{DownarowiczLacroix1996ToeplitzSpectrum}).
Furthermore, efforts have also been made to generalize the notions of Toeplitz sequences and Toeplitz flows (for which the acting group is the group $\Z$ of the integers) to more general classes of group actions.
In the generalized setting, the respective objects are usually referred to as "Toeplitz arrays" and "Toeplitz subshifts".
In \cite{Downarowicz1997RoyalCouple}, the concept of a Toeplitz $\Z^2$-array was first introduced.
A further generalization towards Toeplitz $\Z^d$-arrays and subshifts is presented in \cite{Cortez2006Toeplitz}, and later in \cite{CortezPetite2008Odometers} the authors extend these notions to subshifts over residually finite groups. \\
In this work, we study Toeplitz $\Z^d$-subshifts in the context of topological entropy and unique ergodicity.
More specifically, it is our goal to prove that for any possible given value of entropy there exists a uniquely ergodic Toeplitz $\Z^d$-subshift admitting this entropy.
Moreover, for this subshift we are able to determine its maximal equicontinuous factor.
Altogether, this amounts to the following:
\begin{namedtheorem}[Theorem A]\label{thm: Theorem A}
Let $k\geq 2$ and consider the alphabet $\Sigma_k = \set{0,\ldots,k-1}$.
Then, for every $0<h<\log k$ there exists a Toeplitz-$\Z^d$-array $x\in (\Sigma_k)^{\Z^d}$ such that the Toeplitz subshift $(\ol{O_{\Z^d}}(x),\Z^d)$ is strictly ergodic and $\htop(\ol{O_{\Z^d}}(x),\Z^d)= h$.
Moreover, the maximal equicontinuous factor of $(\ol{O_{\Z^d}}(x),\Z^d)$ is given by the universal $\Z^d$-odometer.
\end{namedtheorem}
In previous works by various authors, one can find related contributions:
\begin{itemize}
    \item In \cite{Cortez2006Toeplitz} an example of a strictly ergodic Toeplitz $\Z^2$-subshift with positive entropy is given.
    \item Krieger proves in \cite{Krieger2007Toeplitz} that any prescribed value can be obtained as entropy of a Toeplitz $G$-subshift, where $G$ is a residually finite amenable group.
    However, he does not obtain control over the number of invariant measures, so that in particular these examples are not necessarily uniquely ergodic.
    \item In \cite{LackaStraszak2018QuasiUniform}, amongst a wealth of other results, the authors reprove Krieger's Theorem using more conceptual techniques.
    \item A version of Krieger's Theorem for sofic entropy is proved in \cite{Kucharski2020Sofic}.
\end{itemize}
It should be noted that these contributions fall short of giving an answer to our goal in the sense that whenever full control over the value of entropy is obtained, control over unique ergodicity is lost, and vice versa.
Moreover, there is no precise mention on the maximal equicontinuous factor of these systems.\\
We should also point out that the statement of Theorem A for d=1 can be considered well-known.
In the literature there are much stronger results, of which the statement in dimension one is a direct corollary.
One such result deals with possible entropy functions on Choquet simplices and is shown in \cite{DownarowiczSerafin2003PossibleEntropy}.
However, the proof already uses the à priori existence of a strictly ergodic dimension 1 Toeplitz flow with prescribed entropy, which is stated as folklore knowledge.
Another strengthening in the one-dimensional case can be found in \cite{Hoynes2016Toeplitz}.
It establishes full control on entropy, the simplex of invariant measures and the maximal equicontinuous factor.
For the proof, the author uses results and refines the techniques of \cite{GjerdeJohansen2000BratteliVershik, Sugisaki2001Toeplitz, Sugisaki2007OrbitEquiv, Sugisaki2011CantorSystems} using properly ordered Bratteli diagrams.\\
An attempt to adapt the methods from above references to the present case of $d>1$ seems fairly hard.
The main obstactle of the transition from the one-dimensional situation is that the linear order of $\Z$ is lost when going to higher dimensions.
For example, many arguments for $d=1$ rely on the well-orderdness of $\N_0$ or a natural notion of concatenation of words.
Such concepts are hard to extend fully to higher dimensions.\\
Here we follow a different and more direct strategy, which also aims to circumvent these issues.
We use an old construction of Grillenberger, which provides a strictly ergodic (but not necessarily Toeplitz) $\Z$-subshift with prescribed entropy \cite{Grillenberger1973GivenEntropy}.
This approach can be modified in order to make the resulting subshift Toeplitz, which therefore already provides a direct and elementary proof of the result in dimension 1.
We then implement this construction in higher dimensions to obtain the proof of Theorem A.
On the way, we also prove a sufficient criterion for unique ergodicity of minimal $G$-subshifts, where $G$ is a residually finite abelian group.
The main advantage of this criterion is that one only needs to count appearances of words along subgroups $\Gamma_t$ of $G$.
To formulate the result, let $(\Gamma_n)_{n\in\N}$ be a decreasing sequence of finite index subgroups of $G$ with trivial intersection and let $(D_n)_{n\in\N}$ be a suitable sequence of fundamental domains of $G/\Gamma_n$.
\begin{namedtheorem}[Theorem B]
    Let $x\in \Sigma^G$ be a $G$-array such that $\ol{O_G}(x)$ is minimal.
    If for every $t\in\N$ and $B\in W_{D_t}(x)$ there exists a real number $\nu(B)$ such that $\mathrm{ap}(B,C)$ converges to $\nu(B)$ uniformly in $C\in W_{D_s}(x)$ as $s\to\infty$, then the subshift $(\ol{O_G}(x),G)$ is strictly ergodic.\\
    Here, $W_{D_t}(x) = \setv{B\in \Sigma^{D_t}}{B \text{ occurs in } x \text{ at some } \gamma_t \in \Gamma_t}$ and $\mathrm{ap}(B,C)$ is the frequency of $B$ appearing in $C$ along elements of $\Gamma_t$.
\end{namedtheorem}
For $G=\Z^d$ we show the converse implication of Theorem B, which yields a direct generalization of \cite[Thm.\ 1.1]{IwanikLacroix1994NonRegular}.\\

\textit{Structure of the article.}
In Section~\ref{sec: prelims}, we recall basic preliminary notions.
We start Section~\ref{sec: Toeplitz} by introducing the concepts of odometers and Toeplitz subshifts.
Then, we discuss a way to decompose a residually finite, abelian group.
This will be of use for the final construction.
In Section~\ref{subsection: unique ergod and entropy} we prove Theorem B by adapting the methods in \cite{IwanikLacroix1994NonRegular} to the setting of residually finite abelian groups.
With similar techniques, we provide a way to compute the topological entropy of a minimal subshift.
Then, we restrict to $\Z^d$ in order to show the converse direction of Theorem B.
Section~\ref{subsection: fundamental domains} contains a version of \cite[Lem.\ 1.9]{Grillenberger1973GivenEntropy}, which allows us to shrink the underlying alphabet and rescale the entropy.
The remainder of Section 4 is dedicated to the proof of Theorem A.
Finally, we conclude with a brief discussion on immediate consequences of Theorem A in Section~\ref{subsection: consequences}.\\[3pt]
{\bfseries Acknowledgements.} I would like to express my gratitude to Tobias Jäger for the supervision and Lino Haupt for the fruitful discussions and helpful comments during the making of this work.
I am also grateful to the referee whose comments helped to improve this article.
Lastly, I want to thank the Graduate Academy of Friedrich Schiller University Jena for supporting my research project on Toeplitz flows through the Honours Programme, which eventually led to this paper.
The author was supported by grant 552775134 of the German Research Council (DFG).

\section{Preliminaries}\label{sec: prelims}
\subsection{Basic notions}
In what follows, we will exclusively work with countable, discrete, abelian groups.
Hence, for simplicity, all definitions will be phrased in this context.
We will write the operation of an abelian group additively and denote the neutral element by $0$.
We call a pair $(X,G)$ a {\bf topological dynamical system (tds)} if $X$ is a compact Hausdorff space, $G$ is a countable, discrete, abelian group and $G$ acts continuously on $X$.
The latter means that there exists a continuous map $\Phi: G\times X \to X$ such that we have $\Phi(0,x) = x$ and $\Phi(g+h,x) = \Phi(g,\Phi(h,x))$ for all $g,h\in G$ and $x\in X$.
Since most of the time the action $\Phi$ will be clear from the context, we will simply write $gx$ for $\Phi(g,x)$.
We denote the {\bf $G$-orbit} of $x$ by $O_G(x) = \setv{gx}{g\in G}$ and the {\bf orbit closure} by $\ol{O_G}(x)$.
We call the tds $(X,G)$ {\bf minimal} if every orbit is dense in $X$.
\\
It is well known that any abelian group $G$ is {\bf amenable}, that is every tds $(X,G)$ admits at least one $G$-invariant Borel probability measure.
In what follows, a "measure" will always mean a Borel probability measure.
Equivalently, $G$ is amenable iff $G$ admits a (left) {\bf F\o lner sequence} (c.f.\ \cite[Thm.\ 7.3]{Pier1984Amenable}), that is a sequence $(F_n)_{n\in\N}$ of finite subsets of $G$ such that for every $g\in G$ we have
$$ \lim_{n\to\infty} \frac{\#[(g+F_n)\symdiff F_n]}{\# F_n} = 0.$$
Since $G$ is countable and discrete, it is a straightforward consequence of \cite[Appendix (3.K)]{Tempelman1992ErgodicTheorems} that $(F_n)_{n\in\N}$ is a F\o lner sequence iff for every non-empty, finite $K\tm G$ we have
$$\lim_{n\to\infty} \frac{\# \partial_K(F_n)}{\#F_n} = 0,$$
where $\partial_K(F_n) = \setv{g\in G}{-K + g \text{ intersects both } F_n \text{ and } G\setminus F_n}$.
We call $\partial_K(F_n)$ the {\bf $K$-boundary} of $F_n$.
One also calls $(F_n)_{n\in\N}$ a {\bf Van Hove sequence} in this case (see \cite[\S 2.4]{Hauser2021RelativeEntropy} for further reading).
The tds $(X,G)$ is {\bf uniquely ergodic} if it admits exactly one invariant measure, and {\bf strictly ergodic} if $(X,G)$ is also minimal.\\
Let $(X,G)$ and $(Y,G)$ be tds.
If there exists a continuous surjective map $f: X\to Y$ such that for every $g\in G$ and $x\in X$ we have $f(gx) = gf(x)$, then we call $(Y,G)$ a {\bf factor} of $(X,G)$ and conversely $(X,G)$ an {\bf extension} of $(Y,G)$.
The map $f$ is called a {\bf factor map} in this case.
If the set of points in $Y$ with exactly one preimage under $f$ is residual in $Y$, we call $f$ an {\bf almost 1-1 factor map} and $(X,G)$ an {\bf almost 1-1 extension} of $(Y,G)$.
It is well-known that in the case of $(Y,G)$ minimal, $f$ is almost 1-1 iff there exists at least one $y\in Y$ with exactly one preimage under $f$.
When $f$ is bijective, we call $(X,G)$ and $(Y,G)$ {\bf topologically conjugate} and $f$ a {\bf conjugacy}.\\
The following well-known characterization of unique ergodicity is a straightforward generalization of the case for $\Z$-actions (c.f.\ \cite[pp.\ 124--125]{Oxtoby1952ErgodicSets}).
We let $C(X)$ be the set of continuous functions $\varphi: X \to \C$.
\begin{thm}\label{thm: chara. unique ergodicity}
    Let $(F_n)_{n\in\N}$ be a F\o lner sequence of $G$, let $(X,G)$ be a tds and $x\in X$. Then, $(\ol{O_G}(x),G)$ is uniquely ergodic iff for all $\varphi \in C(X)$ the averages $\frac{1}{\# F_n}\sum_{h\in F_n+ g}\varphi(hx)$ converge uniformly in $g\in G$ as $n\to \infty$.
\end{thm}

\subsection{Subshifts}
Let $\Sigma$ be some finite alphabet.
The set $\Sigma^G$ of all maps $x: G\to \Sigma$ endowed with the product topology becomes a compact Hausdorff space.
We call the elements of $\Sigma^G$ {\bf arrays}.
We define an action $G\times \Sigma^G \to \Sigma^G, \: (g,x)\mapsto gx$ via
$$ (gx)(h) = x(h+g) \text{ for } h\in G.$$
This yields a tds $(\Sigma^G,G)$, which we call the {\bf full shift} on $\Sigma$.
If $X\tm \Sigma^G$ is compact and shift invariant, we call $(X,G)$ with the restricted action a {\bf subshift}.\\
Let $(F_n)_{n\in\N}$ be a F\o lner (hence also Van Hove) sequence of $G$ and fix a subshift $(X,G)$.
For a finite subset $M$ of $G$ we call a map $B: M\to \Sigma$ an {\bf $M$-block} or just a block.
We shall write $[B]_M = \setc{x\in X}{x\vert_M = B}$ and call such a set a {\bf cylinder set}.
\begin{rem}\label{rem: chara. unique ergodicity subshift}
    Cylinder sets are closed and open, so that their indicator functions are continuous.
    Furthermore, it is an easy consequence of the Stone-Weierstraß Theorem that the linear span of all indicator functions of cylinder sets lies dense in $C(X)$ w.r.t.\ uniform convergence.
    Thus, with Thm.\ \ref{thm: chara. unique ergodicity} it follows that for $x\in \Sigma^G$ the subshift $(\ol{O_G}(x),G)$ is uniquely ergodic iff for all finite $M\tm G$ and blocks $A\in \Sigma^M$ the averages $\frac{1}{\# F_n} \sum_{h\in F_n+g} \1_{[A]_M}(hx)$ converge uniformly in $g\in G$ as $n\to\infty$.
\end{rem}
We denote by $\theta_M(X)=\#\setv{B\in \Sigma^M}{X \cap [B]_M\neq \emptyset}$ the number of $M$-blocks that {\bf occur} in $X$.
With this, we define the {\bf topological entropy of} $(X,G)$ as
\begin{equation} \label{eq: entropy subshift}
\htop(X,G) = \lim_{n\to\infty} \frac{\log \theta_{F_n}(X)}{\# F_n}.
\end{equation}
It follows from the Ornstein-Weiss Lemma (\cite[Thm.\ 1.1]{Krieger2010OrnsteinWeiss}) that this limit exists and is independent of the choice of the Van Hove sequence $(F_n)_{n\in\N}$.
It should be remarked, that the usual approach is to define topological entropy for general tds instead of just subshifts.
Then, one can prove the above equality \eqref{eq: entropy subshift} as a proposition.
See \cite{Hauser2021RelativeEntropy} and \cite{Hauser2022NoteDelone} for further reading.
For the sake of simplicity, we take \eqref{eq: entropy subshift} as the definition instead.

\section{Toeplitz arrays and Toeplitz subshifts}\label{sec: Toeplitz}
Most of the following notions can be found in \cite{CortezPetite2008Odometers} (also for non-abelian groups).
\subsection{$G$-odometers}
Let $G$ be an infinite, countable, discrete, abelian group.
We say that $G$ is {\bf residually finite} if there exists a decreasing sequence $(\Gamma_n)_{n\in\N}$ of finite index (normal) subgroups of $G$ such that $\bigcap_{n\in\mathbb{N}}\Gamma_n=\{0\}$.

Let $(\Gamma_n)_{n\in\mathbb{N}}$ be a decreasing sequence of finite index subgroups of $G$.
The {\bf $G$-odometer} associated to $(\Gamma_n)_{n\in\mathbb{N}}$ is defined as the inverse limit
\begin{align*}
\overleftarrow{G}&:=\varprojlim_n(G/\Gamma_n,\phi_n)
    =\setv{(\ol{g_n})_{n\in\mathbb{N}}\in\prod_{n\in\mathbb{N}}G/\Gamma_n}{\forall n\in\N:\: \phi_n(\ol{g_{n+1}})=\ol{g_n}},
\end{align*}
where $\phi_n: G/\Gamma_{n+1}\to G/\Gamma_n$ is the natural map induced by the inclusion $\Gamma_{n+1}\tm \Gamma_n$.\\
Moreover, with the product topology and pointwise addition $\overleftarrow{G}$ clearly becomes a compact abelian group.
Similarly, $G$ acts on $\overleftarrow{G}$ by pointwise addition.
We use the term $G$-odometer for both the set $\overleftarrow{G}$ and the tds $(\overleftarrow{G},G)$.
This system is  equicontinuous and minimal, hence it is uniquely ergodic with the unique invariant measure being the Haar measure of the group $\overleftarrow{G}$.
We call $\overleftarrow{G}$ {\bf periodic} if some finite index subgroup of $G$ acts trivially on $\overleftarrow{G}$, or equivalently, if the the sequence $(\Gamma_n)_{n\in\N}$ stabilises.
If we additionally assume that $\bigcap_{n\in\mathbb{N}}\Gamma_n=\{0\}$, it can be easily seen that $\tau: G \to \overleftarrow{G},\: g\mapsto (g+\Gamma_n)_{n\in\N}$ is an injective, continuous group homomorphism with dense image.
Hence, we can identify $G$ as a dense subgroup of $\overleftarrow{G}$.
From now on, we assume that $G$ is a residually finite group.

\begin{lem}[{\cite[Lem.\ 2]{CortezPetite2008Odometers}}]\label{lem: factor odometer}
    Let $\overleftarrow{G}_j=\varprojlim_n(G/\Gamma_n^{(j)},\phi^{(j)}_n)$ be two $G$-odometers for $j\in\set{1,2}$.
    There exists a factor map $\pi:\overleftarrow{G}_1\to\overleftarrow{G}_2$ such that $\pi\left((0+\Gamma_n^{(1)})_{n\in\mathbb{N}}\right)=(0+\Gamma_n^{(2)})_{n\in\mathbb{N}}$ iff for every $n\in\mathbb{N}$ there exists $k_n\in\mathbb{N}$ such that $\Gamma_{k_n}^{(1)}\subseteq \Gamma_n^{(2)}$.
    In this case, $\pi$ is also a group homomorphism.
\end{lem}

\begin{rem}
    The above lemma motivates the following definition:
    We call a $G$-odometer $\overleftarrow{G}$ a {\bf universal odometer} if every $G$-odometer is a factor of $\overleftarrow{G}$.
    From the previous lemma it can be seen that whenever a universal odometer exists, it is unique up to isomorphism of topological groups.
    Moreover, $\overleftarrow{G}$ is a universal odometer iff it is a metrizable profinite completion of $G$.
    It can thus be easily seen that every finitely generated group (in particular $\Z^d$) admits a universal odometer.
    Furthermore, the universal $\Z^d$-odometer can be identified with the direct product of $d$ copies of the profinite completion of $\Z$.
\end{rem}

In the case of $G=\Z^d$ one can make some additional observations regarding odometers.
See \cite{Cortez2006Toeplitz} for an extensive discussion.
We call a sequence $(P_n)_{n\in\N}$ of matrices $P_n\in \Z^{d\times d}$ a {\bf scale} if for each $n\in\N$ we he have $\det(P_n)\neq 0$ and there exist $Q_n\in\Z^{d\times d}$ such that $P_{n+1}=P_n\cdot Q_n$.
In this case, the matrix $Q_n=P_n^{-1}\cdot P_{n+1}$ is clearly unique.
We call $(Q_n)_{n\in\N}$ the {\bf increment} of the scale $(P_n)_{n\in\N}$.

\begin{lem}\label{lem: odo. given by scale}
    Let $G=\Z^d$.
    \begin{enumerate}[label=(\roman*)]
        \item For every decreasing sequence $(\Gamma_n)_{n\in\N}$ of finite index subgroups of $\Z^d$ there exists a scale $(P_n)_{n\in\N}$ such that $\Gamma_n = P_n \Z^d$.
        \item Let $(P_n)_{n\in\N}$ be a scale with increment $(Q_n)_{n\in\N}$.
        Suppose that every $P_n$ is a diagonal matrix.
        Then, every $Q_n$ is also diagonal and we have $\bigcap_{n\in\N} P_n \Z^d = \set{0}$ iff $\min_{i\in\set{1,\ldots,d}}\abs{(P_n)_{i,i}} \xrightarrow{n\to\infty}\infty$, where $(P_n)_{i,i}$ denotes the entry of $P_n$ in the $i$-th row and $i$-th column.
    \end{enumerate}
\end{lem}
\begin{proof}
    (i) It follows from the classification of finitely generated abelian groups that each $\Gamma_n$ is isomorphic to $\Z^d$.
    Hence, we can define $P_n$ by writing in each column of $P_n$ one vector of a $\Z$-basis $\mathcal{B}_n$ of $\Gamma_n$.
    Moreover, since $\Gamma_{n+1}\tm \Gamma_n$, we can write any element of $\mathcal{B}_{n+1}$ as a $\Z$-linear combination of $\mathcal{B}_n$.
    Thus, $(P_n)_{n\in\N}$ is a scale.\\
    (ii) This is easy.
\end{proof}
In view of (i) in the above lemma, we call $(P_n)_{n\in\N}$ and $(\Gamma_n)_{n\in\N}$ {\bf associated} if $\Gamma_n = P_n \Z^d$ for all $n\in\N$.
It should be noted that for given $(\Gamma_n)_{n\in\N}$ the associated scale $(P_n)_{n\in\N}$ is not unique.
\begin{cor}\label{cor: universal odo. on Z^d}
    Let $G=\Z^d$ and let $(P_n)_{n\in\N}$ be a scale of diagonal matrices such that for each $M\in \N$ and $i\in \set{1,\ldots,d}$ there exists $n\in\N$ such that $M$ divides $(P_n)_{i,i}$.
    Let $(\Gamma_n)_{n\in\N}$ be the associated sequence of subgroups and $\overleftarrow{G} := \varprojlim_{n}(\Z^d/\Gamma_n,\phi_n)$ the corresponding odometer.
    Then, $\overleftarrow{G}$ is the universal $\Z^d$-odometer.
\end{cor}

\subsection{Toeplitz subshifts}
Let $\Sigma$ be a finite alphabet.
We say that $x\in\Sigma^G$ is a \textbf{Toeplitz array} if for every $g\in G$ there exists a subgroup $\Gamma$ of $G$ with finite index such that $x(g+\gamma)=x(g)$ for each $\gamma\in \Gamma$.
The subshift $(\ol{O_G}(x),G)$ is called a \textbf{Toeplitz subshift}.
It is well-known that every Toeplitz subshift is minimal (c.f.\ \cite[Prop.\ 5]{CortezPetite2008Odometers}).
In \cite[Thm.\ 2.2]{Williams1984ToeplitzFlows} it is shown that the maximal equicontinuous factor of a Toeplitz-$\Z$-subshift is an odometer and conversely it is well-known that every symbolic, minimal almost 1-1 extension of a $\Z$-odometer is a Toeplitz subshift (c.f.\ \cite[Prop.\ 1]{Markley1975Substitution} or \cite[Thm.\ 6]{DownarowiczLacroix1998Almost1-1}).
It is proved in \cite{CortezPetite2008Odometers} that the same holds in the generalized setting.\\
Let $x\in \Sigma^G$. We define for finite index subgroup $\Gamma$ and $\alpha\in\Sigma$
\begin{align*}
    \Per(x,\Gamma,\alpha)&=\setv{g\in G}{x(g+\gamma)=\alpha \text{ for all }\gamma\in \Gamma},\\
    \Per(x,\Gamma)&=\bigcup_{\alpha\in \Sigma}\Per(x,\Gamma,\alpha) = \setv{g\in G}{x(g+\gamma) = x(g) \text{ for all } \gamma\in \Gamma}.
\end{align*}
Thus, as $G$ is countable, we see that $x$ is a Toeplitz array iff there exists a decreasing sequence $(\Gamma_n)_{n\in\mathbb{N}}$ of finite index subgroups of $G$ such that $G=\bigcup_{n\in\mathbb{N}}\Per(x,\Gamma_n)$.
Assume that $x\in\Sigma^G$ is a Toeplitz array.
We call $\Gamma$ an {\bf essential period} of $x$, if $\Per(x,\Gamma)$ is non-empty and the implication
$$(\forall \alpha\in\Sigma: \:\Per(x,\Gamma,\alpha)\subseteq \Per(x,\Gamma,\alpha)-h) \implies h\in \Gamma$$
holds for all $h\in G$.
A decreasing sequence $(\Gamma_n)_{n\in\N}$ of finite index subgroups of $G$ is called a \textbf{period structure} of $x$ if $G=\bigcup_{n\in\mathbb{N}}\Per(x,\Gamma_n)$ and every $\Gamma_n$ is an essential period.
It is shown in \cite[Cor.\ 6]{CortezPetite2008Odometers} that every Toeplitz array has a period structure.

\begin{prop}[{\cite[Prop.\ 7]{CortezPetite2008Odometers}}]\label{prop: MEF of Toeplitz}
    Let $x\in\Sigma^G$ be a Toeplitz array with period structure $(\Gamma_n)_{n\in\mathbb{N}}$ and let $\overleftarrow{G}$ the $G$-odometer associated to $(\Gamma_n)_{n\in\mathbb{N}}$.
    Then, $(\overleftarrow{G},G)$ is the maximal equicontinuous factor (MEF) of the Toeplitz subshift $(\ol{O_G}(x),G)$ via the the factor map $\pi:\ol{O_G}(x) \to\overleftarrow{G}$ defined by
    \begin{align*}
        \pi(y)=(g_n+\Gamma_n)_{n\in \N}, \text{ iff } \Per(y,\Gamma_n,\alpha) = \Per(x,\Gamma_n,\alpha)-g_n \text{ for all } \alpha \in \Sigma \text{ and } n\in\N.
    \end{align*}
    Moreover, the map $\pi$ satisfies $\# \pi^{-1}(\pi(y)) = 1$ iff $y\in \ol{O_G}(x)$ is a Toeplitz array.
    In particular, $\pi$ is almost 1-1.
\end{prop}

\subsection{Decomposing a residually finite, abelian group}\label{subsection: theta maps}
For a subgroup $\Gamma$ of $G$ a we call a set $D\tm G$ a {\bf fundamental domain} of $G/\Gamma$ if each coset in $G/\Gamma$ has exactly one representative $d\in D$, or in other words, if the canonical projection $D\to G/\Gamma, \: d\mapsto d+\Gamma$ is a bijection.
\begin{lem}[{\cite[Lem.\ 5]{CortezPetite2014Invariant}}]\label{lem: ex. of fundamental domains amenable}
    Let $(\Gamma_n)_{n\in\N}$ be a decreasing sequence of subgroups of $G$ with finite index such that $\bigcap_{n\in\N} \Gamma_n = \set{0}$.
    Then, there exists a strictly increasing sequence of integers $(n_k)_{k\in\N}$ and a sequence $(D_k)_{k\in\N}$ of subsets of $G$ such that for all $k\in \N$ we have
    \begin{enumerate}[label=(\roman*)]
        \item $G = \bigcup_{k\in\N} D_k$.
        \item $0\in D_k \tm D_{k+1}$.
        \item $D_k$ is a fundamental domain of $G/\Gamma_{n_k}$.
        \item $D_l = \bigcup_{\gamma\in D_l \cap \Gamma_{n_k}} (\gamma +D_k)$ for $l>k$.
        \item $(D_k)_{k\in\N}$ is a F\o lner sequence in $G$.
    \end{enumerate}
\end{lem}
It should be noted that analogues of this result hold for the case where $G$ is non-abelian or even non-amenable.
Suppose we are given for $(\Gamma_n)_{n\in\N}$ with $\bigcap_{n\in\N} \Gamma_n =\set{0}$, sequences $(n_k)_{k\in\N}$ and $(D_k)_{k\in\N}$ according to Lem.\ \ref{lem: ex. of fundamental domains amenable}.
By going over to a subsequence, we can assume w.l.o.g.\ that $n_k=k$ and $\Gamma_{k} \neq \Gamma_{k-1}$ for all $k\in\N$, i.e.\ $D_k\cap \Gamma_{k-1} \neq \set{0}$.
We also let $\Gamma_0 = G$.
\begin{defn}\label{def: theta maps}
Let $g\in G$.
By Lem.\ \ref{lem: ex. of fundamental domains amenable} (i), there exists $n:= \min\setv{i\in\N}{g\in D_i}\in\N$ and by (iv), we can write $g$ uniquely as $g = \gamma_{n} + d_{n-1}$ with $\gamma_n\in D_n \cap \Gamma_{n-1}$ and $d_{n-1} \in D_{n-1}$.
Now we can apply the same argument to $d_{n-1}$ to obtain a decomposition $d_{n-1} = \gamma_{n-1} + d_{n-2}$.
Iterating this yields a unique decomposition
$$g = \sum_{i=1}^n \gamma_i \text{ with }\gamma_i\in D_i \cap \Gamma_{i-1}.$$
We define $\theta_i(g) = \gamma_i$ for $i\in \set{1,\ldots,n}$ and $\theta_i(g) = 0$ for $i>n$.
Thus, we obtain maps $\theta_i : G \to D_i\cap \Gamma_{i-1}$ for $i\in \N$.
\end{defn}

\begin{ex}\label{ex: theta maps}
    Let us give an illustration for the previous definition.
    We consider subgroups $\Gamma_n= P_n\Z^2$ ($n\in \set{1,\ldots,4}$) of $G=\Z^2$ with
    \begin{align*}
         &P_1  = \diag(2,2) \quad (\text{green dots}),\\
         &P_2 = \diag(4,6) \quad (\text{blue rings}),\\
         &P_3 = \diag(8,12) \quad(\text{orange rings}),\\
         &P_4 = \diag(24,24),
    \end{align*}
    and fundamental domains, which are coloured according to their respective subgroups,
    \begin{align*}
        & D_1 = \setv{h=(h_1,h_2)\in \Z^2}{0\leq h_1 \leq 1 \:\text{ and } -1\leq h_2 \leq 0}\\
        & D_2 = \setv{h=(h_1,h_2)\in \Z^2}{-2\leq h_1 \leq 1 \:\text{ and } -3\leq h_2 \leq 2}\\
        & D_3 = \setv{h=(h_1,h_2)\in \Z^2}{-2\leq h_1 \leq 5 \:\text{ and } -9\leq h_2 \leq 2}\\
        & D_4 = \setv{h=(h_1,h_2)\in \Z^2}{-10\leq h_1 \leq 13 \:\text{ and } -9\leq h_2 \leq 14}.
    \end{align*}
    We compute for $g = (10,5)^T$ that
    \begin{align*}
        \theta_4(g) = \begin{pmatrix}
            8\\
            12
        \end{pmatrix}
        ,\quad
        \theta_3(g) = \begin{pmatrix}
            4\\
            -6
        \end{pmatrix}
        , \quad
        \theta_2(g) = \begin{pmatrix}
            -2\\
            0
        \end{pmatrix}
        , \quad
        \theta_1(g) = \begin{pmatrix}
            0\\
            -1
        \end{pmatrix}.
    \end{align*}
    A visualization of this computation can be found below in Figure \ref{fig: theta maps}.
    \begin{figure}
    \centering
    \includegraphics[width=1\linewidth]{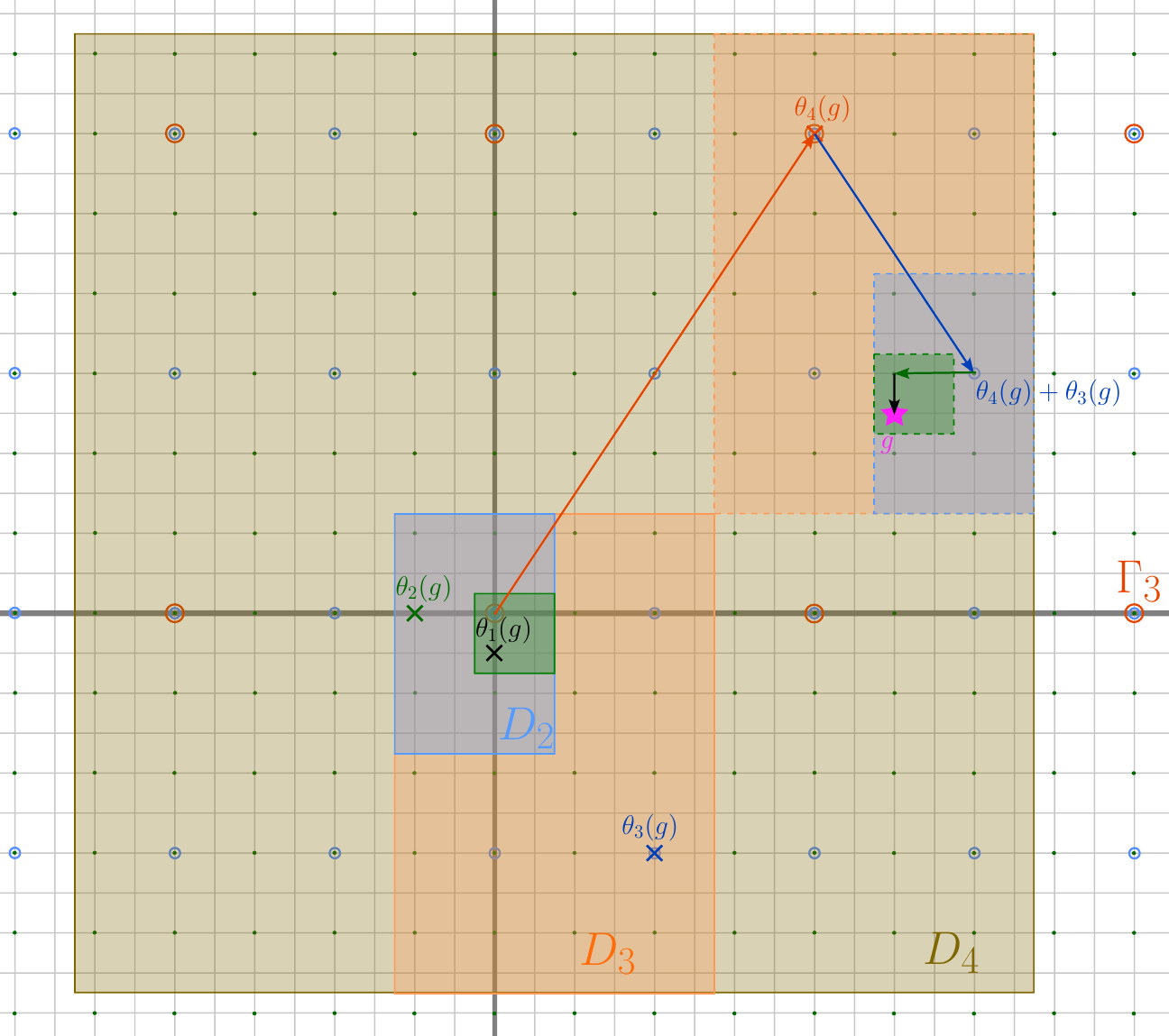}
    \caption{\it A visualization of Example \ref{ex: theta maps}.
    Shifted copies of the fundamental domains are marked by dashed borders.}
    \label{fig: theta maps}
    \end{figure}
\end{ex}

The following properties are immediate consequences from the definition:
\begin{lem}\label{lem: properties theta maps}
Let $g\in G$.
\begin{enumerate}[label=(\roman*)]
    \item If $g\in D_m$, then $g = \sum_{i=1}^m \theta_i(g)$.
    \item We have $g\in \Gamma_n$ iff $\theta_i(g) = 0$ for all $i\in\set{1,\ldots,n}$.
    \item If $\gamma\in \Gamma_n$, then $\theta_n(g+\gamma) = \theta_n(g)$.
    \item If $\theta_i(g)+\theta_i(h) \in D_i$ for all $i\in\set{1,\ldots,n}$, then $\theta_n(g+h) = \theta_n(g)+\theta_n(h)$.
\end{enumerate}
\end{lem}

The following lemma will be used later to show that the constructed periods are essential.
Let $x\in \Sigma^G$ be an arbitrary array.
\begin{lem}\label{lem: property * implies essential period}
    Let $h\in G$ and let $\Gamma$ be a finite index subgroup of $G$.
    Assume that there exists some finite set $I \tm \N$ such that
    $$ \Per(x,\Gamma) = \setv{g\in G}{\theta_i(g)=0 \text{ for some } i\in I}.$$
    If $\Per(x,\Gamma,\alpha)\tm \Per(x,\Gamma,\alpha)-h$ holds for all $\alpha\in \Sigma$, then $\theta_i(h) =0$ for all $i\in I$.
\end{lem}
\begin{proof}
    Assume for a contradiction that $\theta_j(h)\neq 0$ for some $j\in I$.
    We can define an element $g:= \sum_{i\in I} \gamma_i \in G$ as follows:
    Whenever $\theta_i(h)\neq 0$, we let $\gamma_i=0$.
    Conversely, if $\theta_i(h) = 0$, then we choose an arbitrary $\gamma_i\in (D_i\cap\Gamma_{i-1})\setminus\set{0}$.
    With this we have on one hand $g\in \Per(x,\Gamma,\alpha)$ with $\alpha = x(g)$, since $\theta_j(g) = \gamma_j = 0$.
    However, it follows by Lem.\ \ref{lem: properties theta maps} (iv) for $i\in I$ that
    $$ \theta_i(g+h) = \theta_i(h) \text{ if } \theta_i(h)\neq 0 \text{ and } \theta_i(g+h) = \gamma_i \neq 0 \text{ if } \theta_i(h) = 0.$$
    Thus, $g+h\notin \Per(x,\Gamma)$ which contradicts $\Per(x,\Gamma,\alpha) \tm \Per(x,\Gamma,\alpha)-h$.
\end{proof}

\subsection{Unique ergodicity and entropy for Toeplitz subshifts}\label{subsection: unique ergod and entropy}
\begin{figure}
        \centering
        \includegraphics[width=.7\linewidth]{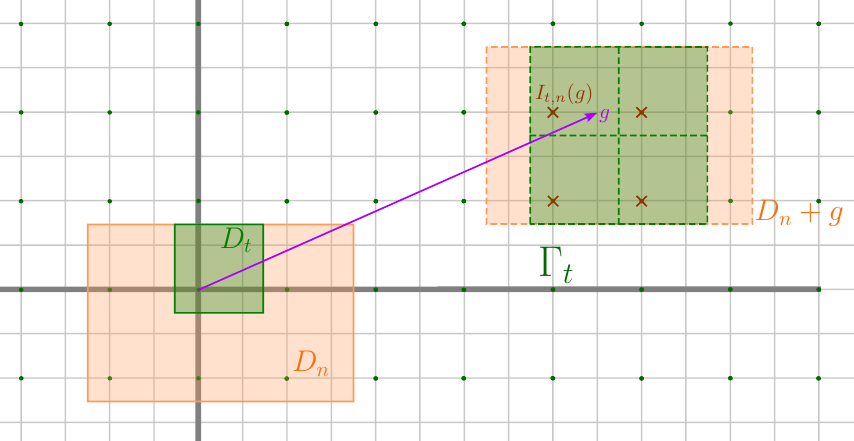}
        \caption{\it The decomposition of $D_n+g$ into $\biguplus_{\gamma\in I_{t,n}(g)} (D_t+\gamma)$ and $R_{t,n}(g)$.
        The green dots are the elements of the subgroup $\Gamma_t$ and the brown crosses are the points in the set $I_{t,n}(g)$.}
        \label{fig: interior and rest}
    \end{figure}

In what follows, we will use the ideas of \cite{IwanikLacroix1994NonRegular} to obtain means to ensure unique ergodicity while also having control over the topological entropy.
Let $(\Gamma_n)_{n\in\N}$ be a decreasing sequence of finite index subgroups of $G$ with $\bigcap_{n\in\N} \Gamma_n = \set{0}$.
Let $x\in \Sigma^G$ be an array such that the subshift $(\ol{O_G}(x),G)$ is minimal.
(At this point, it is not necessary that $x$ is Toeplitz.)
Furthermore, let $(n_k)_{n\in\N}$ and $(D_k)_{k\in\N}$ be given according to Lem.\ \ref{lem: ex. of fundamental domains amenable}.
We assume w.l.o.g.\ that $n_k = k$.
We call a $D_t$-block $B: D_t\to \Sigma$ a {\bf $D_t$-symbol} of $x$, if there exists $\gamma \in \Gamma_t$ such that $B(d) = x(d+\gamma)$ for all $d\in D_t$.
We denote the set of all $D_t$-symbols of $x$ by $W_{D_t}(x)$.\\
Let $s,t\in \N$ with $s\geq t$ and let $B\in W_{D_t}(x)$, $C\in W_{D_s}(x)$.
We define the {\bf frequency of $B$ appearing in $C$} as
$$ \mathrm{ap}(B,C) = \frac{1}{\#D_s\cap \Gamma_t} \cdot \#\setv{\gamma\in D_s\cap \Gamma_t}{\forall d\in D_t: C(d+\gamma)=B(d)}.$$
Note that we do not count all the occurrences of $B$ in $C$, but only occurrences along points in $\Gamma_t$.
The following is Theorem B.
\begin{thm}\label{thm: sufficient crit. unique ergod. Toeplitz}
    If for every $t\in\N$ and $B\in W_{D_t}(x)$ there exists a real number $\nu(B)$ such that $\mathrm{ap}(B,C)$ converges to $\nu(B)$ uniformly in $C\in W_{D_s}(x)$ as $s\to\infty$, then the subshift $(\ol{O_G}(x),G)$ is strictly ergodic.
\end{thm}
Here, we mean with "$\mathrm{ap}(B,C)$ converges to $\nu(B)$ uniformly in $C\in W_{D_s}(x)$ as $s\to\infty$" that for every $\varepsilon>0$ there exists $s_0\in \N$ such that for every $s\geq s_0$ and $C\in W_{D_s}(x)$ we have $\abs{\mathrm{ap}(B,C)-\nu(B)}<\varepsilon$.
Also note that since $\ap(B,C)$ does not count all occurences of $B$ in $C$, the result is not a direct consequence of Rem.\ \ref{rem: chara. unique ergodicity subshift}\\[3pt]
    The proof will be divided below into several lemmas.
    Let $M\tm G$ be finite and $A\in \Sigma^M$ be a block.
    For this proof we will use the notation
    \begin{align*}
    &x\vert_{M+h} = A :\iff \forall m\in M: x(m+h) = A(m) \text{ and}\\
    &B\vert_{M+h} = A :\iff \forall m\in M: B(m+h) = A(m),
    \end{align*}
    whenever $B\in \Sigma^{M^\prime}$ with $M+h \tm M^\prime$.
    In view of Rem.\ \ref{rem: chara. unique ergodicity subshift}, it suffices to show that
    \begin{align*} F(A,g,n) &:= \frac{1}{\# D_n} \sum_{h\in D_n+g} \1_{[A]_M}(hx)
        = \frac{1}{\# D_n} \#\{h\in D_n + g : x\vert_{M+h}=A\}
    \end{align*}
    converges uniformly in $g\in G$ as $n\to \infty$.
    We assume w.l.o.g.\ that $0\in M$.
    Let $\varepsilon>0$ and choose $t$ such that $\frac{\# \partial_{-M}(D_t)}{\# D_t}<\frac{\varepsilon}{4}$.
    For a $D_t$-symbol $B$ we let $N(A,B)$ be the number of occurrences of $A$ in $B$, i.e.\ $N(A,B)= \# \setv{h\in G}{M+h \tm D_t \text{ and } B\vert_{M+h}=A}$.
    Given $n\in\N$ and $t$ as above, write $D_n + g = \biguplus_{\gamma\in I_{t,n}(g)} (D_t+\gamma) \uplus R_{t,n}(g)$,
    where we define
    $$ I_{t,n}(g) = \setv{\gamma\in \Gamma_t}{D_t + \gamma \tm D_n + g} \quad \text{and}\quad R_{t,n}(g) = (D_n + g)\setminus \biguplus_{\gamma\in I_{t,n}(g)} (D_t + \gamma).$$
    See Figure \ref{fig: interior and rest} for an illustration in the case $G=\Z^2$.

    \begin{lem}\label{lem: limit R_{t,n}(g)}
        We have
        $$ \lim_{n\to\infty} \frac{\# R_{t,n}(g)}{\# D_n} = 0 \text{ and } \lim_{n\to\infty} \frac{\# I_{t,n}(g)}{\# D_n} = \frac{1}{\# D_t}$$
        and both convergences are uniform in $g\in G$.
    \end{lem}
    \begin{proof}
    We claim that $\setv{h\in D_n+g}{(D_t-D_t)+h\tm D_n +g}\tm \biguplus_{\gamma\in I_{t,n}(g)}(D_t+\gamma)$.
    Let $h\in D_n+g$ such that $(D_t-D_t)+h \tm D_n+g$ and let $\gamma\in \Gamma_t$ be such that $h\in D_t+\gamma$.
    Then, $\gamma-h \in -D_t$ and hence
    $$ D_t+\gamma = D_t + \gamma - h + h \tm (D_t-D_t)+ h \tm D_n+g.$$
    Thus, $\gamma\in I_{t,n}(g)$, proving the claim.\\
    With the claim, it is now easy to see that $R_{t,n}(g)\tm \partial_{D_t-D_t}(D_n+g)$.
    Since $(D_n)_{n\in\N}$ is a Van Hove sequence and $D_t-D_t$ is finite, it follows that
    $$\frac{\# R_{t,n}(g)}{\# D_n} \leq \frac{\# \partial_{D_t-D_t}(D_n+g)}{\#D_n} = \frac{\# \partial_{D_t-D_t}(D_n)}{\# D_n} \xrightarrow{n\to\infty} 0,$$
    proving $\frac{\# R_{t,n}(g)}{\# D_n} \xrightarrow{n\to\infty} 0$ uniformly in $g$.
    Moreover, since $\# D_n = \#(D_n+g) = \# D_t \cdot \# I_{t,n}(g) + \# R_{t,n}(g)$, the second part follows easily.
    \end{proof}

    As a next step we claim that
    \begin{lem}
    Under the assumptions of Thm.\ \ref{thm: sufficient crit. unique ergod. Toeplitz} we have
    $$ \Psi(t,n,B,g):= \frac{1}{\# I_{t,n}(g)} \#\setc{\gamma\in I_{t,n}(g)}{x\vert_{D_t+\gamma}=B} \xrightarrow{n\to\infty} \nu(B)$$
    for all $B\in W_{D_t}(x)$ and the convergence is uniform in $g\in G$.
    \end{lem}
    \begin{proof}
    Let $\delta>0$, $B\in W_{D_t}(x)$ and choose $s\geq t$ such that $\abs{\mathrm{ap}(B,C)-\nu(B)} < \delta/4$ for all $C\in W_{D_s}(x)$.
    Since $D_s = \biguplus_{\gamma_t\in D_s\cap \Gamma_t}(D_t+\gamma_t)$, we obtain the decomposition
    \begin{equation} \label{eq: D_t D_s D_n}
    D_n+g = \kl\biguplus_{\gamma_s \in I_{s,n}(g)} \biguplus_{\gamma_t \in D_s \cap \Gamma_t} (D_t +\gamma_t + \gamma_s)\kr \uplus R_{s,n}(g).
    \end{equation}
    It is furthermore not hard to see for $\gamma\in I_{t,n}(g)$ that we either have $\gamma \in (D_s\cap \Gamma_t)+ I_{s,n}(g)$ or $D_t+\gamma \tm R_{s,n}(g)$.
    (For an illustration see Figure \ref{fig: D_t D_s D_n} below.)
    This yields $\abs{\Psi(t,n,B,g)-\nu(B)} \leq U_n(g) + V_n(g)$ with
    \begin{align*}
    U_n(g) &= \abs{\frac{1}{\# I_{t,n}(g)} \sum_{\gamma_s\in I_{s,n}(g)} \#\setc{\gamma_t\in D_s\cap \Gamma_t}{x\vert_{D_t+\gamma_t+\gamma_s}=B}-\nu(B)}
    \quad \text{ and }\\
    V_n(g) &= \frac{1}{\# I_{t,n}(g)} \#\setv{\gamma_t\in \Gamma_t}{D_t+\gamma_t\tm R_{s,n}(g)\text{ and } x\vert_{D_t+\gamma_t}=B}.
    \end{align*}
    Observe that
    $$V_n(g) \leq \frac{\# R_{s,n}(g)}{\# I_{t,n}(g)\cdot \# D_t}
    = \frac{\# R_{s,n}(g)}{\# D_n}\cdot \frac{\# D_n}{\# I_{t,n}(g)\cdot \# D_t} < \frac{\delta}{4}$$
    for large enough $n$ by the previous lemma.
    Moreover, we have
    $$ \# \setc{\gamma_t \in D_s\cap \Gamma_t}{x\vert_{D_t+\gamma_t+\gamma_s}=B} = \#(D_s \cap \Gamma_t) \cdot \mathrm{ap}(B,C_{\gamma_s}),$$
    where $C_{\gamma_s}: D_s \to \Sigma,\: d\mapsto x(d+\gamma_s)$ is a $D_s$-symbol.
    Thus,
    \begin{align*}
    U_n(g) &= \abs{\kl\frac{\#(D_s\cap \Gamma_t)}{\# I_{t,n}(g)} \sum_{\gamma_s\in I_{s,n}(g)} \mathrm{ap}(B,C_{\gamma_s})\kr-\nu(B)} \\
        &\leq \frac{\# (D_s\cap\Gamma_t)}{\# I_{t,n}(g)} \sum_{\gamma_s\in I_{s,n}(g)} \abs{\mathrm{ap}(B,C_{\gamma_s})-\nu(B)} \\
        &\hspace{2cm}+ \abs{\frac{\# (D_s\cap \Gamma_t)\cdot \# I_{s,n}(g)}{\# I_{t,n}(g)}\nu(B) - \nu(B)} \\
        &\leq \frac{\# (D_s\cap \Gamma_t)\cdot \# I_{s,n}(g)}{\# I_{t,n}(g)}\cdot \frac{\delta}{4} + \nu(B)\abs{\frac{\# (D_s\cap \Gamma_t)\cdot \# I_{s,n}(g)}{\# I_{t,n}(g)} - 1} \\
        &< \frac{3\delta}{4}.
    \end{align*}
    for $n$ large enough, as
    \[ \frac{\# (D_s\cap \Gamma_t)\cdot \# I_{s,n}(g)}{\#I_{t,n}(g)} = \frac{\# D_s\cdot\# I_{s,n}(g)}{\# D_n} \cdot \frac{\#D_n}{\# D_t \cdot \# I_{t,n}(g)} \xrightarrow{n\to\infty} 1. \qedhere\]
    \end{proof}

    \begin{figure}
    \centering
    \includegraphics[width=.7\linewidth]{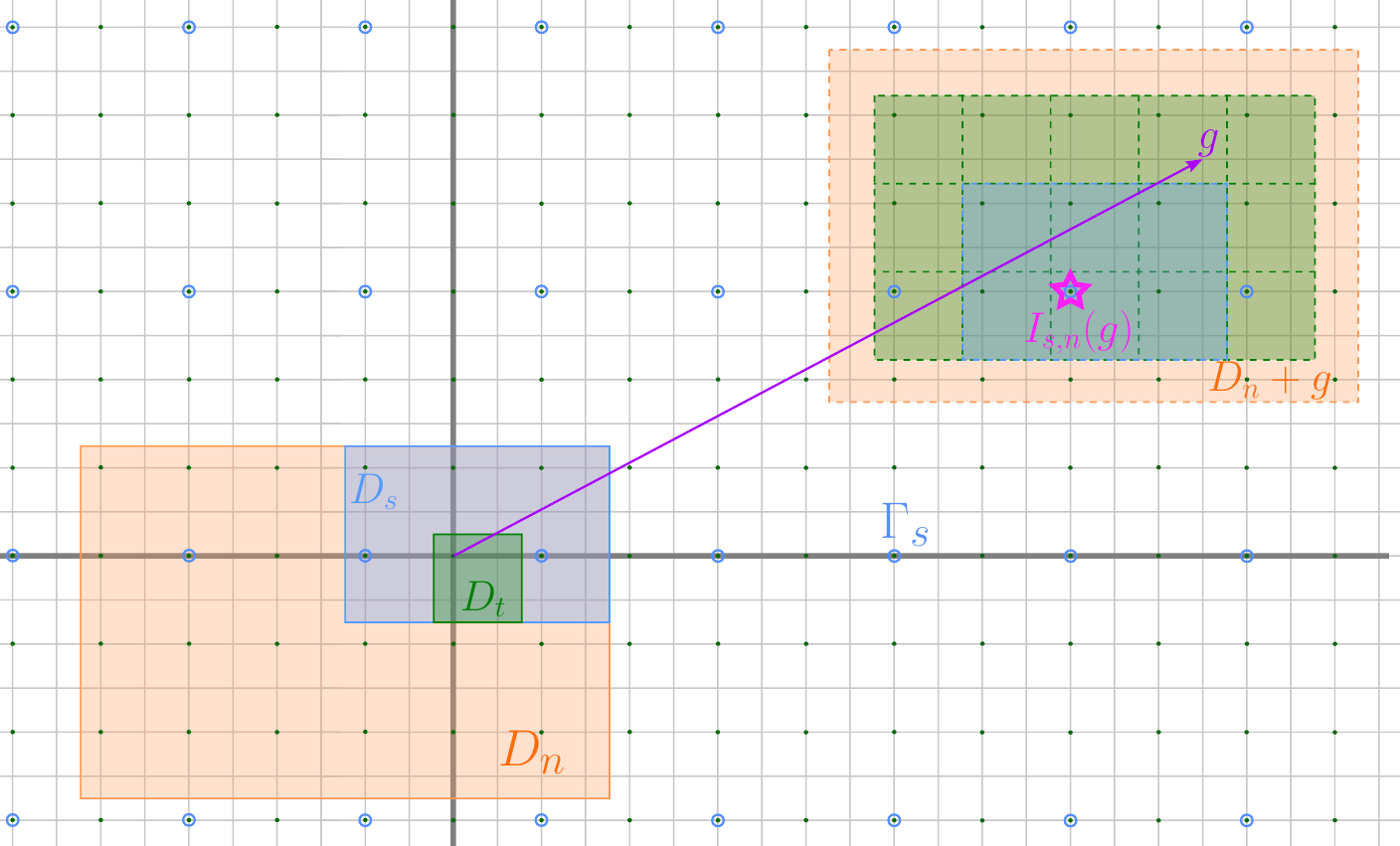}
    \caption{\it A visualization of the decomposition in \eqref{eq: D_t D_s D_n}. In this particular setting, the set $I_{s,n}(g)$ consists of a single point marked by the pink star.
    The dark green, dashed area corresponds to the double union in \eqref{eq: D_t D_s D_n}, whereas the light green, dashed squares are the sets $D_t+\gamma \tm R_{s,n}(g)\:(\gamma\in \Gamma_t)$.}
    \label{fig: D_t D_s D_n}
    \end{figure}

    We continue with the proof of Thm.\ \ref{thm: sufficient crit. unique ergod. Toeplitz}.
    Choose $\delta>0$ such that
    $$\frac{\delta}{\# D_t} \sum_{B\in W_{D_t}(x)} N(A,B)< \frac{\varepsilon}{4}.$$
    By the above lemma, the difference between $\Psi( t,n,B^\prime,g)$ and $\nu(B^\prime)$ is less than $\delta$ for all $B^\prime\in W_{D_t}(x)$ and $g\in G$ if $n$ is large enough.
    Let $h\in\setc{h^{\prime\prime}\in D_n+g}{x\vert_{M+h^{\prime\prime}}= A}$.
    We can write $h = h^\prime + \gamma$ uniquely with $\gamma\in \Gamma_t$ and $h^\prime\in D_t$.
    Let $B_\gamma$ be the $D_t$-symbol given by $B_\gamma(d) =x(d+\gamma)$ for $d\in D_t$.
    Assume that $\gamma\in I_{t,n}(g)$ and $M+h^\prime\tm D_t$.
    It follows for all $m\in M$ that  $A(m)=x(m+h)=x(m+h^\prime+\gamma) =B_\gamma(m+h^\prime)$.
    Hence, $h^\prime \in \setv{h^{\prime\prime}\in G}{M+h^{\prime\prime}\tm D_t\text{ and } (B_\gamma)\vert_{M+h^{\prime\prime}}=A}$.
    This shows for all $g\in G$
    \begin{align*}
        F(A,g,n) &=  \frac{1}{\# D_n} \#\setc{h\in D_n + g}{x\vert_{M+h}=A}
        \geq \frac{1}{\# D_n}\sum_{\gamma\in I_{t,n}(g)} N(A,B_\gamma)\\
        &= \frac{1}{\# D_n}\sum_{B\in W_{D_t}(x)} \#\setc{\gamma\in I_{t,n}(g)}{x\vert_{D_t+\gamma}=B} \cdot N(A,B)\\
        &\geq \frac{\# I_{t,n}(g)}{\# D_n} \sum_{B\in W_{D_t}(x)}(\nu(B)-\delta)N(A,B)\\
        &= \frac{\# I_{t,n}(g)}{\# D_n} \sum_{B\in W_{D_t}(x)} \nu(B) N(A,B) - \frac{\# I_{t,n}(g)}{\# D_n}\cdot \delta \sum_{B\in W_{D_t}(x)} N(A,B)\\
        &\geq \frac{1}{\# D_t} \sum_{B\in W_{D_t}(x)} \nu(B)N(A,B) - \varepsilon,
    \end{align*}
    if $n$ is chosen large enough.\\
    In order to obtain an upper bound on $F(A,g,n)$, note that the occurrences $h = h^\prime + \gamma$ ($h^\prime \in D_t, \: \gamma\in \Gamma_t$) which we have not counted in the above way must {\bf not} satisfy $\gamma\in I_{t,n}(g)$ and $M+h^\prime\tm D_t$.
    Hence, if $\gamma\in I_{t,n}(g)$, then $M+h^\prime \cap (G\setminus D_t)\neq \emptyset$.
    Since $h^\prime \in D_t$ and $0\in M$, this implies $h^\prime \in \partial_{-M}(D_t)$.
    Therefore, there are at most $\# I_{t,n}(g) \cdot \# \partial_{-M}(D_t)$ occurrences of this form.
    On the other hand, if $M+h^\prime\tm D_t$, then $\gamma\notin I_{t,n}(g)$.
    Since $h\in D_n+g$, it follows that $h\in R_{t,n}(g)$.
    Thus, there can be at most $\# R_{t,n}(g)$ occurrences of this form.
    Altogether we obtain for all $g\in G$
    \begin{align*}
        F(A,g,n) &\leq \frac{\# I_{t,n}(g)}{\# D_n} \sum_{B\in W_{D_t}(x)}(\nu(B)+\delta)N(A,B) \\
        &\hspace{1cm}+\frac{\# I_{t,n}(g)}{\# D_n} \# \partial_{-M}(D_t) + \frac{\# R_{t,n}(g)}{\# D_n}\\
        &< \frac{\# I_{t,n}(g)}{\# D_n} \sum_{B\in W_{D_t}(x)}\nu(B)N(A,B) + \frac{\# I_{t,n}(g)}{\#D_n} \cdot \delta \sum_{B\in W_{D_t}(x)}N(A,B) + \frac{\varepsilon}{2}\\
        &< \frac{1}{\# D_t} \sum_{B\in W_{D_t}(x)} \nu(B)N(A,B)+ \varepsilon,
    \end{align*}
    if $n$ is large enough.
    Hence, we obtain for $m,n\in \N$ large enough and all $g\in G$ that $\abs{F(A,g,m)-F(A,g,n)} < 2\varepsilon$.
    This is a uniform Cauchy condition which yields the desired uniform convergence.
    This finishes the proof of Thm.\ \ref{thm: sufficient crit. unique ergod. Toeplitz}.\qed \\

As a particular case, we obtain the following.

\begin{cor}\label{cor: sufficient crit. unique ergod. Toeplitz}
    If $\mathrm{ap}(B,C) = \mathrm{ap}(B,C^\prime)$ for any $B\in W_{D_t}(x)$ and $C,C^\prime \in W_{D_{t+1}}(x)$ ($t=1,2,\ldots$), then the subshift $(\ol{O_G}(x),G)$ is strictly ergodic.
\end{cor}

Now we turn our attention to entropy.
In a similar way as above, a simple estimate shows that instead of dealing with all blocks over the alphabet $\Sigma$ (c.f.\ \eqref{eq: entropy subshift}), it suffices to restrict oneself to the $D_t$-symbols.
\begin{prop}\label{prop: entropy Toeplitz}
    We have
    $$ \htop(\ol{O_G},G) = \lim_{n\to\infty} \frac{1}{\# D_n} \log \#W_{D_n}(x).$$
\end{prop}
\begin{proof}
    Recall that we denote for $n\in \N$
    $$ \theta_{D_n}(x):=\theta_{D_n}(\ol{O_G}(x)) = \# \setv{B\in \Sigma^{D_n}}{\exists g\in G: \forall d\in D_n: x(d+g) = B(d)}.$$
    Thus, we clearly have $\theta_{D_n}(x) \geq \# W_{D_n}(x)$ which yields
    $$ \htop(\ol{O_G}(x),G) \geq \limsup_{n\to \infty} \frac{1}{\# D_n} \log\# W_{D_n}(x).$$
    For the converse inequality, we use a similar argument as in the proof of Thm.\ \ref{thm: sufficient crit. unique ergod. Toeplitz}.
    Let $t,s\in \N$ with $s\geq t$.
    Then, using the partition
    $D_s+g = \biguplus_{\gamma\in I_{t,s}(g)}(D_t+\gamma)\uplus R_{t,s}(g)$
    for $g\in G$, one sees quickly that
    $$ \theta_{D_s}(x) \leq \sup_{g\in G}\kl\#\Sigma^{\# R_{t,s}(g)} \cdot (\# W_{D_t}(x))^{\# I_{t,s}(g)}\kr < \infty.$$
    Due to Lem.\ \ref{lem: limit R_{t,n}(g)}, it follows for every $\varepsilon>0$ that
    \begin{align*}
        \frac{\log \theta_{D_s}(x)}{\#D_s}&\leq\sup_{g\in G}\kl \frac{\#R_{t,s}(g)\log \# \Sigma}{\# D_s}+\frac{\# I_{t,s}(g)\log \# W_{D_t}(x)}{\# D_s}\kr \\
        &\leq \frac{1}{\# D_t}\log \# W_{D_t}(x)+\varepsilon
    \end{align*}
    for large enough $s$.
    Hence, $\htop(\ol{O_G}(x),G) \leq \frac{1}{\# D_t} \log \# W_{D_t}(x)$.
    This implies
    \[ \htop(\ol{O_G}(x),G)\leq \liminf_{n\to\infty} \frac{1}{\#D_n} \log \# W_{D_n}(x).\qedhere\]
\end{proof}

It should be noted that even though Thm.\ \ref{thm: sufficient crit. unique ergod. Toeplitz} and Prop.\ \ref{prop: entropy Toeplitz} are valid for arbitrary minimal subshifts, they are most applicable in a situation where there is a canonical choice for the sequence $(\Gamma_n)_{n\in\N}$.
Most importantly, if $x$ is a Toeplitz array, we can choose $(\Gamma_n)_{n\in\N}$ as a period structure of $x$.
In this situation, it turns out that under additional assumptions, it is possible to prove the converse of Thm.\ \ref{thm: sufficient crit. unique ergod. Toeplitz}.
To that end, we now assume that $G=\Z^d$ and that $x$ is a Toeplitz array with period structure $(\Gamma_n)_{n\in\N}$.
Let $(P_n)_{n\in\N}$ be an associated scale and $\overleftarrow{G}$ be the associated odometer.
Observe that since subsequences of period structures are period structures, we can assume w.l.o.g.\ that the sequence $(n_k)_{k\in\N}$ given by Lem.\ \ref{lem: ex. of fundamental domains amenable} satisfies $n_k = k$.
In this case, it is possible to express the Toeplitz subshift generated by $x$ as a certain skew product system.
See \cite{IwanikLacroix1994NonRegular} for the construction in the case of $G=\Z$.\\
For $t\in \N$, we define maps $\pi_t: \overline{O_G}(x) \to D_t$ via
$$ \pi_t(y) = d \text{ iff } \Per(y,\Gamma_t,\alpha) = \Per(x,\Gamma_t,\alpha)-d \text{ for all } \alpha \in \Sigma.$$
With this we have $\pi(y)= (\pi_t(y)+\Gamma_t)_{t\in\N}$, where $\pi: \overline{O_G}(x)\to \overleftarrow{G}$ is the almost 1-1 factor map from Prop.\ \ref{prop: MEF of Toeplitz}.
In particular, each $\pi_t$ is continuous.
Fix $t\in \N$, let $y\in \overline{O_G}(x)$ and $\gamma\in \Gamma_t$.
We write $W_t(y,\gamma)$ for the $D_t$-block which is defined by
$$ W_t(y,\gamma)(d)=y(d+\gamma-\pi_t(y)) \text{ for } d\in D_t.$$
The following is an easy consequence of writing $y = \lim_{n\in\N} g_n x$ with $g_n \in G$.
\begin{lem}
    The block $W_t(y,\gamma)$ is a $D_t$-symbol of $x$.
\end{lem}

We can now define $y^{(t)}\in W_{D_t}(x)^{G}$ by setting
$$ y^{(t)}(g) = W_t(y, P_t g) \text{ for } g\in G.$$
Of course, the mapping $y\mapsto y^{(t)}$ is continuous.
We let $\psi_t: G \to D_t$ be the canonical projection defined for $g\in G$ by
\begin{equation}\label{eq: def psi maps}
    \psi_t(g) = d \iff d\in D_t \text{ and } g+\Gamma_t = d+\Gamma_t.
\end{equation}
This is well defined since $D_t$ is a fundamental domain of $G/\Gamma_t$.
Define the map
$$ \epsilon_t: G\times G \to G,\: (g,h)\mapsto P_t^{-1}(-\psi_t(g+h)+g+\psi_t(h) ).$$
When applied to $D_t\times D_t$, one can interpret $\epsilon_t$ as some sort of carry over.
\begin{rem}
    Suppose that $G=\Z$ and let $\Gamma_t = P_t \Z$ with some integer $P_t\geq 2$.
    We naturally choose $D_t = \{0,1,\ldots, P_t-1\}$ as the fundamental domain of $G/\Gamma_t$.
    We aim to compute $\epsilon_t(c,d)$ for $c,d\in D_t$.
    Assume first that $d< P_t-c$.
    Then,
    $$\epsilon_t(c,d) = \frac{1}{P_t}(-\psi_t(c+d) + c + \psi_t(d)) = \frac{1}{P_t}(-(c+d) +c + d) = 0.$$
    We compute furthermore for $d\geq P_t-c$, i.e.\ $d = P_t-c + d^\prime$ with $d^\prime \in \set{0,\ldots, c-1}$:
    $$ \epsilon_t(c,d) = \frac{1}{P_t}(-\psi_t(P_t+d^\prime) + c + \psi_t(d)) = \frac{1}{P_t}(-d^\prime+c+d) = \frac{P_t}{P_t} = 1.$$
    We have shown that
    $$ \epsilon_t(c,d) = \begin{cases}
        0   &, d < P_t-c\\
        1   &, \text{else}
    \end{cases}.$$
    Note that in \cite{IwanikLacroix1994NonRegular}, the authors define a related map $\varepsilon_t: \Z/P_t\Z \to \{0,1\}$ that satisfies $\varepsilon_t(d+P_t\Z)=\epsilon_t(1,d)$ and plays a similar role in their construction in the case $G=\Z$.
\end{rem}
We define further
$$ \delta_t: G\times G/\Gamma_t \times \overline{O_G}(x^{(t)})\to \overline{O_G}(x^{(t)}),\: [\delta_t(g,h+\Gamma_t,y)](g^\prime) = y(g^\prime + \epsilon_t(g,h)).$$
In other words, the triplet $(g,h+\Gamma_t,y)$ gets mapped to the sequence $y$ shifted by $\epsilon_t(g,h)$.
Note that $\epsilon_t(g,h)$ is independent of the representative $h$ of the coset $h+\Gamma_t$.
This allows us to define a skew product tds $(G/\Gamma_t \times \overline{O_G}(x^{(t)}), G)$, where the $G$-action is defined for $g,h\in G$ and $y\in \overline{O_G}(x^{(t)})$ via
$$  (g,(h+\Gamma_t, y))\mapsto  (g+h+\Gamma_t, \delta_t(g,h+\Gamma_t,y)).$$
It is not hard to see that $G$ acts indeed by homeomorphisms on $G/\Gamma_t \times \overline{O_G}(x^{(t)})$.
With this, we aim to ensure that any "carry over" from the addition in the first coordinate will be taken into account via an appropriate shift in the second coordinate.
\begin{rem}
    Above skew product action is a direct generalization of the corresponding homeomorphism defined in \cite[p.\ 193]{IwanikLacroix1994NonRegular} (denoted as $\widetilde{S_t}$ therein).
\end{rem}
\begin{lem}\label{lem: Toeplitz conjugate to skew prod.}
\hspace{.5cm}
    \begin{enumerate}[label=(\roman*)]
        \item If $y\in \overline{O_G}(x)$, then $y^{(t)}\in \overline{O_G}(x^{(t)})$.
        \item The map $\Phi_t: \overline{O_G}(x)\to (G/\Gamma_t)\times \overline{O_G}(x^{(t)}),\: y\mapsto (\pi_t(y)+\Gamma_t, y^{(t)})$ is a homeomorphism.
        \item The Toeplitz subshift $(\overline{O_G}(x),G)$ is topologically conjugate to the skew product system $(G/\Gamma_t\times \overline{O_G}(x^{(t)}),G)$ via the map $\Phi_t$.
        \item If the Toeplitz subshift $(\overline{O_G}(x),G)$ is uniquely ergodic, then so is the subshift $(\ol{O_G}(x^{(t)}),G)$.
    \end{enumerate}
\end{lem}
\begin{proof}
    (i) This follows again by writing $y = \lim_{n\to\infty} g_n x$ and using continuity of $\pi_t$.
    (ii) Continuity of $\Phi_t$ is clear.
    For injectivity, assume that $\pi_t(y)+\Gamma_t =\pi_t(z)+\Gamma_t$ and $y^{(t)}=z^{(t)}$.
    The first equality implies $\pi_t(y)=\pi_t(z)=: d_t$, and combined with the second equality this yields for all $\gamma\in\Gamma_t$ and $d\in D_t$ that  $y(d+\gamma-d_t) = z(d+\gamma-d_t)$.
    Since $D_t+\Gamma_t = G$, it follows that $y=z$.
    For surjectivity, let $\xi = \lim_{n\to\infty} g_n x^{(t)}\in \ol{O_G}(x^{(t)})$ and let $h+\Gamma_t \in G/\Gamma_t$.
    By compactness, we can assume that $((\psi_t(h)+P_t g_n)x)_{n\in\N}$ converges to some $y\in \ol{O_G}(x)$, where $\psi_t$ is the projection defined by \eqref{eq: def psi maps}.
    It is not hard to check that $\Phi_t(y) = (h+\Gamma_t, \xi)$.\\
    (iii) Let $g\in G$ and $y \in \ol{O_G}(x)$.
    Since we have $\pi_t(gy) = \psi_t(g+\pi_t(y))$ (so that $\pi_t(gy)+\Gamma_t = (g+\pi_t(y))+\Gamma_t$), it remains to show that $(gy)^{(t)} = \delta_t(g,\pi_t(y)+\Gamma_t,y^{(t)})$.
    For this let $h\in G$ and $d\in D_t$.
    It follows that
    \begin{align*}
        [(gy)^{(t)}(h)](d) &= W_t(gy,P_th)(d)
        = (gy)(d+P_t h -\pi_t(gy)) \\
        &= y(d+P_t h - \psi_t(g+\pi_t(y)) + g)\\
        &= y(d+P_t h +\underbrace{(-\psi_t(g+\pi_t(y)) + g + \pi_t(y))}_{=P_t \epsilon_t(g,\pi_t(y))} - \pi_t(y))\\
        &= y(d+P_t(h+\epsilon_t(g,\pi_t(y))) - \pi_t(y))\\
        &= W_t(y,P_t(h+\epsilon_t(g,\pi_t(y))))(d)
        = [y^{(t)}(h+\epsilon_t(g,\pi_t(y)))](d)
    \end{align*}
    This shows the desired equality
    $$
        (gy)^{(t)}(h) =y^{(t)}(h+\epsilon_t(g,\pi_t(y))) = [\delta_t(g,\pi_t(y)+\Gamma_t,y^{(t)})](h).
    $$
    (iv) This follows immediately from (iii).
\end{proof}

Now we are in position to prove the converse to Thm.\ \ref{thm: sufficient crit. unique ergod. Toeplitz} in the case $G=\Z^d$.
\begin{thm}\label{thm: necessary crit. unique ergod. Toeplitz}
    If the Toeplitz subshift $(\ol{O_G}(x),G)$ is strictly ergodic, then for every $t\in\N$ and $B\in W_{D_t}(x)$ there exists a real number $\nu(B)$ such that $\mathrm{ap}(B,C)$ converges to $\nu(B)$ uniformly in $C\in W_{D_s}(x)$ as $s\to\infty$.
\end{thm}
\begin{proof}
    Let $t\in \N$ and $B\in W_{D_t}(x)$.
    If we define $F_n= P_t^{-1}(D_n\cap \Gamma_t)$, it is not hard to see that $(F_n)_{n\in\N}$ is a F\o lner sequence.
    By Lem.\ \ref{lem: Toeplitz conjugate to skew prod.}(iv), the subshift $(\ol{O_G}(x^{(t)}),G)$ is uniquely ergodic.
    Thus, by Thm.\ \ref{thm: chara. unique ergodicity} there exists some constant $\nu(B)$ such that
    \begin{align*} &\frac{1}{\# F_s}  \#\setv{h\in F_s + g}{x^{(t)}(h) = B} \\
    &\hspace{1cm} = \frac{1}{\#(D_s\cap \Gamma_t)} \#\setv{h\in P_t^{-1}(D_s\cap \Gamma_t)}{x^{(t)}(h+g) = B}
    \end{align*}
    converges uniformly in $g\in G$ to $\nu(B)$ as $s\to \infty$.
    (Recall that we treat $B \in W_{D_t}(x)$ as a letter of the array $x^{(t)}$.)
    For every $s\geq t$ and $C\in W_{D_s}(x)$ there exists some $\gamma_C\in \Gamma_s \tm \Gamma_t = P_t\Z^d$ such that $C(d) = x(d+\gamma_C)$ for $d\in D_s$.
    Set $g_C = P_t^{-1}\gamma_C$.
    Now we see that for $h\in F_s = P_t^{-1}(D_s\cap \Gamma_t)$,
    \begin{align*}
        x^{(t)}(h+g_C) = B &\iff \forall d\in D_t: x(d+P_t h + P_t g_C) = B(d) \\
        &\iff \forall d\in D_t: x(d+ P_t h+\gamma_C) = B(d)\\
        &\iff \forall d\in D_t: C(d+P_t h) = B(d).
    \end{align*}
    This shows that
    $$ \frac{1}{\# F_s} \#\setv{h\in F_s + g_C}{x^{(t)}(h) = B} = \mathrm{ap}(B,C),$$
    which finishes the proof.
\end{proof}

\section{The construction}\label{sec: construction}

\subsection{Obtaining a fundamental domain}\label{subsection: fundamental domains}
Recall that Lem.\ \ref{lem: ex. of fundamental domains amenable} yields for a given decreasing sequence $(\Gamma_n)_{n\in\N}$ of finite index subgroups with $\bigcap_{n\in\N}\Gamma_n =\set{0}$, an appropriate family of fundamental domains of a subsequence of the quotient groups $G/\Gamma_n$.
Before, we have ignored the unpleasant restriction that the lemma only provides the result for a subsequence, by going over to said subsequence.
This was possible, since we did not rely on any properties of $\Gamma_n$ and $D_n$ outside the ones given by this lemma.
The main problem why this does not suffice for our construction, is that in order to define blocks $C: D_n\to \Sigma_k$, we need to know exactly the index of $\Gamma_n$ in $G$, i.e.\ the cardinality of $D_n$.
Fortunately, one sees that in the situation of our construction, such a family $(D_n)_{n\in\N}$ for the full sequence $(\Gamma_n)_{n\in\N}$ is not hard to obtain.
We will give a short outline for one way to achieve this.\\
We assume for the entirety of Section 4 that $G=\Z^d$. Recall that for the proof of Theorem A we wish to construct the universal $\Z^d$-odometer as the MEF of our Toeplitz subshift.
We already know by Cor.\ \ref{cor: universal odo. on Z^d}, that the universal odometer can be expressed via a scale $(P_n)_{n\in\N}$ consisting of diagonal matrices.
Moreover, we can assume w.l.o.g.\ that all diagonal entries of $P_n$ are positive.
Another advantage is that in geometric terms, this situation is as simple as possible.
For these reasons, we will restrict our view to this setting for the remainder of this section.\\
Let $(P_n)_{n\in\N}$ be a scale of diagonal matrices $P_n \in \Z^{d\times d}$ with positive diagonal entries.
Let $(Q_n)_{n\in\N}$ be the increment of $(P_n)_{n\in\N}$, let $(\Gamma_n)_{n\in\N}$ be the associated sequence of subgroups and assume that $\bigcap_{n\in\N} \Gamma_n =\set{0}$.
We interpret each matrix $P_n$ as a linear map $\R^d \to \R^d,\: \xi \mapsto P_n\xi$.
Therefore, we can consider the set
\begin{align*}
    P_n([0,1)^d) &= \setv{P_n\xi}{\xi=(\xi_1,\ldots,\xi_d)\in \R^d, \forall i\in \set{1,\ldots,d}: 0\leq \xi_i < 1} \\
        &= \setv{\zeta\in \R^d}{\forall i\in\set{1,\ldots,d}: 0\leq \zeta_i < (P_n)_{i,i}}.
\end{align*}
We now define
$$D_n^\prime = P_n([0,1)^d) \cap \Z^d = \setv{g=(g_1,\ldots,g_d)\in \Z^d}{\forall i \in \set{1,\ldots,d}: 0\leq g_i <(P_n)_{i,i}}.$$
Note that $D_n^\prime$ is just a $d$-dimensional rectangle which sits in the "first quadrant" of the coordinate system.
Moreover, Lem.\ \ref{lem: odo. given by scale} yields that all the sidelengths of $D_n^\prime$ go to infinity as $n\to\infty$.
Therefore, the following is clear from a geometric standpoint.
\begin{lem}\label{lem: cuboids are good fundamental domains}
    We have for all $n\in \N$:
    \begin{enumerate}[label=(\roman*)]
        \item $0\in D_n^\prime \tm D_{n+1}^\prime$.
        \item $D_n^\prime$ is a fundamental domain of $G/\Gamma_n$.
        \item $D_m^\prime = \bigcup_{\gamma\in D_m^\prime \cap \Gamma_n} (\gamma + D_n^\prime)$ for $m>n$.
        \item $(D_n^\prime)_{n\in\N}$ is a F\o lner sequence in $G$.
    \end{enumerate}
\end{lem}

It is now easy to see that shifting each $D_n^\prime$ by an appropriate amount to obtain some set $D_n$ will imply the property $G=\bigcup_{n\in\N} D_n$.
One only needs to be careful to choose the right amount, so that the other properties are preserved.
The following solution to this is very simple and - most importantly - applicable to our later construction.
\begin{lem}\label{lem: shifting cuboids for better fundamental domain}
    Assume that there exists an $N\in \N$ such that we have $(Q_n)_{i,i}\geq 4$ for all $n\geq N$ and $i\in \set{1,\ldots,d}$.
    Define, for $n\in\N$, vectors $r_n,s_n \in \Z^d$ as follows:
    For $i\in \set{1,\ldots,d}$ let $(r_n)_i = \lfloor (Q_n)_{i,i}/4\rfloor$. Moreover, let $s_1=0$ and $s_{n+1} = P_n r_n + s_n$.
    If we set $D_n = D_n^\prime - s_n$ for $n\in\N$, then
    \begin{enumerate}[label=(\roman*)]
            \item $G = \bigcup_{n\in\N} D_n$.
            \item $0\in D_n\tm D_{n+1}$.
            \item $D_n$ is a fundamental domain of $G/\Gamma_n$.
            \item $D_m = \bigcup_{\gamma\in D_m \cap \Gamma_n} (\gamma + D_n)$ for $m>n$.
            \item $(D_n)_{n\in\N}$ is a F\o lner sequence in $G$.
        \end{enumerate}
\end{lem}
\begin{proof}
    (i) It follows by induction on $n\in \N$, that for all $i\in\set{1,\ldots,d}$ we have $(s_n)_i \leq \frac{(P_n)_{i,i}}{2}$.
    Furthermore, we have $(s_n)_i \geq (P_{n-1})_{i,i}$ for $n>N$.
    Since we defined $D_n^\prime = \setv{g\in \Z^d}{\forall i\in \set{1,\ldots,d}: 0\leq g_i < (P_n)_{i,i}}$, above observations yield
    $$ \setv{g\in\Z^d}{\forall i\in\set{1,\ldots,d}: -(P_{n-1})_{i,i}\leq g_i < \frac{(P_n)_{i,i}}{2}} \tm D_n$$
    for all $n>N$.
    Since $(P_n)_{i,i}\xrightarrow{n\to\infty}\infty$, the first assertion follows. \\
    (ii),(iv) Note that $P_n r_n \in D_{n+1}^\prime\cap \Gamma_n$.
    With this it easily follows by induction that $s_n \in D_n^\prime$.
    Therefore, $0\in D_n^\prime - s_n = D_n$.
    Moreover, we also obtain
    $$ D_{n+1} = D_{n+1}^\prime - s_{n+1} = \bigcup_{\gamma\in D_{n+1}^\prime \cap \Gamma_n} (\gamma -P_n r_n+ (D_n^\prime - s_n)).$$
    Since $\gamma\in D_{n+1}^\prime \cap \Gamma_n \iff \gamma-P_nr_n \in D_{n+1}\cap \Gamma_n$, we obtain (iv) and thus also (ii).
    Assertions (iii) and (v) are clear.
\end{proof}
Another tool for proving Theorem A is being able to transform a Toeplitz array by replacing its entries by blocks $B\in \Sigma^D$ with a suitable finite set $D\tm G$.
In the context of $\Z$-subshifts, this procedure is known as a {\bf substitution} (of constant length).
Grillenberger has shown that in the case of an injective constant length substitution, strict ergodicity is preserved and entropy is divided by the length of the substitution (c.f.\ \cite[Lem.\ 1.9]{Grillenberger1973GivenEntropy}).\\
Assume that $(\Gamma_n)_{n\in\N}$ is a period structure of the Toeplitz array $x$ and that the scale $(P_n)_{n\in\N}$ satisfies the conditions of Lem. \ref{lem: shifting cuboids for better fundamental domain}.
Let $P\in\Z^{d\times d}$ be another diagonal matrix with positive diagonal entries.
We also set $\Gamma = P\Z^d$ and let $D = P([0,1)^d)\cap \Z^d$ be a fundamental domain of $G/\Gamma$.
Let $S: \Sigma\to\Sigma^{D}$ be an injective map and $y\in \Sigma^G$.
We define, for $\gamma\in \Gamma$ and $d\in D$,
$$ y^S(d+\gamma) = [S(y(P^{-1}\gamma))](d)$$
This extends $S$ to a map $\Sigma^G\to \Sigma^G,\: y\mapsto y^S$.
Similarly, let $B\in \Sigma^M$ with finite $M\tm G$.
Then we can set for $\gamma\in PM$ and $d\in D$
$$ B^S(d+\gamma) = [S(B(P^{-1}\gamma))](d),$$
which extends $S$ to a map $\Sigma^M\to \Sigma^{D + P M},\: B\mapsto B^S$.
Clearly, these extensions remain injective.
An important observation is that $(P\cdot P_n)_{n\in\N}$ remains a scale with the same increment as $(P_n)_{n\in\N}$.
We let $\Delta_n= P\Gamma_n= (PP_n)\Z^d$ denote the subgroup associated to $P\cdot P_n$.
Moreover, if we set $F_n^\prime = (PP_n)([0,1)^d)\cap \Z^d$ and $F_n = F_n^\prime - Ps_n$ (where $s_n$ is defined as in Lem.\ \ref{lem: shifting cuboids for better fundamental domain}), we see that $F_n^\prime = D + P D_n^\prime$ as well as $F_n = D + PD_n$.
Hence, it easily follows that $\Delta_n$ and $F_n$ satisfy properties (i) - (v) of Lem.\ \ref{lem: shifting cuboids for better fundamental domain}.

\begin{rem}\label{rem: period structure after substitution}
    Even if $(\Gamma_n)_{n\in\N}$ is a period structure of a Toeplitz array $x$, it is in general not true that $(\Delta_n)_{n\in\N}$ is a period structure of $x^S$.
    Examples can be already found in the case of $G=\Z$, where $\Delta_1$ is not an essential period of $x^S$.
    However, in the case of $\bigcap_{n\in\N} \Gamma_n = \set{0}$ it seems difficult to decide whether or not it is true in general that $\varprojlim_{n}(G/\Delta_n,\phi_n)$ is the MEF of $(\ol{O_G}(x^S),G)$.
    For the specific sequence $x$ which we construct later, we will find an elementary argument to prove that $(\Delta_n)_{n\in\N}$ is indeed a period structure of $x^S$.
\end{rem}

\begin{prop}\label{prop: entropy and strict ergod. under substitution}
    If $x$ is a Toeplitz array generating a strictly ergodic Toeplitz subshift, then so is $x^S$.
    Moreover, we have
    $$ \htop(\ol{O_G}(x^S),G) = \frac{1}{\# D} \cdot\htop(\ol{O_G}(x),G).$$
\end{prop}
\begin{proof}
    Observe that if $g\in \Per(x,\Gamma_n)$, then $D+Pg\tm \Per(x^S,\Delta_n)$.
    As $D+P\Z^d=\Z^d$, this shows that $x^S$ is a Toeplitz array.
    Furthermore, it is not hard to see that
    $$ W_{F_n}(x^S) = \setv{B^S}{B\in W_{D_n}(x)},$$
    so that $\# W_{F_n}(x^S) =\# W_{D_n}(x)$.
    Together with $\# F_n = \# D\cdot \# D_n$, it follows that
    \begin{align*}
        \htop(\ol{O_G}(x^S),G) &= \lim_{n\to\infty} \frac{1}{\# F_n} \log \# W_{F_n}(x^S)\\
                &= \lim_{n\to\infty} \frac{1}{\# D}\cdot \frac{1}{\# D_n} \log \# W_{D_n}(x)\\
                &= \frac{1}{\# D}\cdot \htop(\ol{O_G}(x),G).
    \end{align*}
    Concerning strict ergodicity, let $t\in\N$ and $s\geq t$.
    Also let $B^\prime\in W_{F_t}(x^S)$ and $C^\prime\in W_{F_{s}}(x^S)$.
    By what we have shown, there exist $B \in W_{D_t}(x)$ and $C \in W_{D_{s}}(x)$ such that $B^\prime=B^S$ and $C^\prime=C^S$.
    With a similar argument one can see that
    \begin{align*}
        \mathrm{ap}(B^\prime,C^\prime) &:= \frac{1}{\#(F_{s}\cap \Delta_t)} \cdot \# \setv{\delta \in F_{s}\cap \Delta_t}{\forall f\in F_t: C^\prime(f+\delta)=B^\prime(f)}\\
        &= \frac{1}{\#(D_{s}\cap \Gamma_t)}\cdot \#\setv{\gamma\in D_{s}\cap \Gamma_t}{\forall d\in D_t: C(d+\gamma)=B(d)}\\
        &= \mathrm{ap}(B,C).
    \end{align*}
    By strict ergodicity of $(\ol{O_G}(x),G)$ and the assumption that $(\Gamma_n)_{n\in\N}$ is a period structure of $x$, it follows via Thm.\ \ref{thm: necessary crit. unique ergod. Toeplitz} that $\mathrm{ap}(B^\prime,C^\prime)$ converges uniformly in $C^\prime \in W_{F_{s}}(x^S)$ to some real number $\nu(B^\prime)$ as $s \to\infty$.
    Thus, Thm.\ \ref{thm: sufficient crit. unique ergod. Toeplitz} yields strict ergodicity of $(\ol{O_G}(x^S),G)$.
\end{proof}

\subsection{Three sequences}\label{subsection: three sequences}
We shall now construct three sequences $(p^\prime_n)_{n\in\N_0}$, $(q^\prime_n)_{n\in\N_0}$ and  $(\lambda_n)_{n\in\N_0}$ of real numbers.
This construction is similar to the one in the beginning of Section 2 of \cite{Grillenberger1973GivenEntropy}.
It will turn out that $p^\prime_n$ and $q^\prime_n$ will be used for defining diagonal matrices $P_n, Q_n \in \Z^{d\times d}$ and $\lambda_n $ will be related to the topological entropy.\\
For $k\geq 5$ we define $p^\prime_n := p^\prime_n(k)$ and $q^\prime_n:= q^\prime_n(k)$ inductively by
\begin{align*}
    &p^\prime_0 = 1            && q^\prime_0=k \\
    &p^\prime_{n+1} = p^\prime_n q^\prime_n  && q^\prime_{n+1}=(q^\prime_n -1)!
\end{align*}
Now we also let $\lambda_n := \lambda_n(k) := \frac{1}{p^\prime_n} \log(q^\prime_n)$, so that $q^\prime_n = e^{p^\prime_n\lambda_n}$.
We recall that for any $m\in \N$, it follows from the representation of $e^m$ via the exponential series that
\begin{equation}\label{eq: estimate factorial}
    1 < (m!)m^{-m}e^m.
\end{equation}
The following lemma bears lots of similarities with \cite[Lem.\ 2.1]{Grillenberger1973GivenEntropy}.
\begin{lem}\label{lem: properties p_n, lambda_n}
    For any $k\geq 5$ we have the following:
    \begin{enumerate}[label=(\roman*)]
        \item The sequence $(\lambda_n)_{n\in\N_0}$ is monotonically decreasing.
        In particular, the limit $\lambda(k):= \lim_{n\to\infty} \lambda_n = \inf_{n\in\N_0} \lambda_n$ exists.
        \item $\lambda(k)>\log k - 5$.
        \item $\lambda(k)>0$.
    \end{enumerate}
\end{lem}
\begin{proof}
    (i) We use the fact that $m!\leq m^m$ for any positive integer $m$ to obtain
    $$ e^{p^\prime_{n+1} \lambda_{n+1}} = q^\prime_{n+1} < (q^\prime_n)! \leq {q^\prime_n}^{q^\prime_n} = e^{p^\prime_n \lambda_n q^\prime_n} = e^{p^\prime_{n+1} \lambda_n},$$
    which yields monotonicity of $(\lambda_n)_{n\in\N_0}$.\\
    (ii) Plugging $m=q^\prime_j-1$ into \eqref{eq: estimate factorial} and taking the natural logarithm yields
    \begin{align*}0 &< \log(q^\prime_{j+1}) - (q^\prime_j -1) \log(q^\prime_j-1) + q^\prime_j-1 \\
        &= p^\prime_{j+1}\lambda_{j+1} - (q^\prime_j -1) \log(q^\prime_j-1) + q^\prime_j-1.
    \end{align*}
    Using the mean value theorem on the natural logarithm, we can see that
    $$ p^\prime_j\lambda_j - \log(q^\prime_j-1) = \log(e^{p^\prime_j\lambda_j}) - \log(e^{p^\prime_j\lambda_j}-1) \leq \frac{1}{e^{p^\prime_j\lambda_j}-1} \leq 1,$$
    and hence $ p^\prime_j \lambda_j - 1 \leq \log(q^\prime_j-1)$.
    Therefore,
    \begin{align*}
        0   &<p^\prime_{j+1}\lambda_{j+1} - (q^\prime_j-1)(p^\prime_j\lambda_j-1) + q^\prime_j-1
        \leq p^\prime_{j+1}\lambda_{j+1} -p^\prime_{j+1}\lambda_j +p^\prime_j\lambda_j + 2q^\prime_j -2.
    \end{align*}
    Dividing by $p^\prime_{j+1}$ yields
    $$ 0 < \lambda_{j+1} - \lambda_j + \frac{p^\prime_j\lambda_j}{p^\prime_{j+1}} + 2\kl \frac{1}{p^\prime_j} - \frac{1}{p^\prime_{j+1}}\kr.$$
    Since $k\geq 5$, it follows by induction on $j\geq 1$ that $p^\prime_j = p^\prime_j(k) \geq 2$ as well as $p^\prime_j\lambda_j \geq j$.
    Moreover, using $e^{-1/2}\leq 2/3$, we see that
    $$ \frac{p^\prime_j\lambda_j}{p^\prime_{j+1}} = \frac{\lambda_j}{e^{p^\prime_j\lambda_j}} \leq e^{\lambda_j-p^\prime_j\lambda_j}\leq e^{(-1/2)p^\prime_j\lambda_j} \leq \kl\frac{2}{3}\kr^j$$
    for $j\geq 1$.
    Also note that $\frac{p^\prime_0\lambda_0}{p^\prime_1}= \frac{\log k}{k} \leq 1 = \kl\frac{2}{3}\kr^0$, which yields for all $j\in \N_0$
    \begin{equation}\label{eq: estimate lambda_(j+1)-lambda_j}
        0 < \lambda_{j+1}-\lambda_{j} + \kl\frac{2}{3}\kr^j + 2\kl\frac{1}{p^\prime_j}-\frac{1}{p^\prime_{j+1}}\kr.
    \end{equation}
    Summing \eqref{eq: estimate lambda_(j+1)-lambda_j} over $j=0,1,2,\ldots$ now gives us
    $$ 0 < \lambda(k) - \lambda_0 + 3 + \frac{2}{p^\prime_0} = \lambda(k) - \log k + 5.$$
    (iii) By (ii), we deduce that $\lambda(k)>0$ for $k\geq e^5$.
    Moreover, if one sets $\tilde{k} = (k-1)!$, then it easily follows that $q^\prime_j(\tilde{k}) = q^\prime_{j+1}(k)$ and $kp^\prime_j(\tilde{k}) = p^\prime_{j+1}(k)$.
    This implies $\lambda_j(\tilde{k}) = k\cdot\lambda_{j+1}(k)$ and hence
    \begin{equation}\label{eq: recursive formula lambda(k)}
        \lambda(k) = \frac{1}{k} \cdot \lambda(\tilde{k}) = \frac{1}{k}\cdot \lambda((k-1)!).
    \end{equation}
    As $(k-1)! > k$ holds for all $k\geq 5$, iterating \eqref{eq: recursive formula lambda(k)} often enough yields $\lambda(k)>0$.
\end{proof}

\subsection{Proof of Theorem A}
Let us first give a rough outline of the proof, which consists of two major steps.
First, we show that any entropy $h$ can be realized, as long as we additionally assume $k\geq 5$ and $h<\lambda(k)$.
For this, we will first construct sets of blocks which will be the $D_n$-symbols of the constructed array $x$.
All the $D_{n+1}$-blocks are constructed as a "concatenation" or "permutation" of all the previously constructed $D_n$-blocks.
This will yield strict ergodicity.
Moreover, the relation $\# W_{D_{n+1}}(x) \approx \# W_{D_n}(x)!$ will imply an exponential growth rate of the number of $D_n$-symbols and hence positive entropy.
With careful choice of the exact growth (see Lem.\ \ref{lem: properties q_n} below), we will obtain precise control over the entropy.\\
Comparing this intermediate result with Theorem A, we are left with a "gap" of possible entropies $h\in [\lambda(k),\log k)$.
However, Lem.\ \ref{lem: properties p_n, lambda_n} (ii) yields that the length of this gap is always bounded.
In particular, the relative length of this gap compared to the whole interval $(0,\log k)$ becomes arbitrarily small for growing $k$.
Therefore, by choosing some integer $s$ such that $k^s$ is large enough, we can realize the entropy $sh$ with some Toeplitz array over the alphabet $\Sigma_{k^s}$.
Now, with an appropriate substitution $S: \Sigma_{k^s} \to (\Sigma_k)^F$ where $F\tm \Z^d$ with $\# F = s$, it follows by Prop.\ \ref{prop: entropy and strict ergod. under substitution} that $y := x^S$ does the job.
We also note that the character "$0$" is used for the real number zero, a symbol from the alphabet $\Sigma_k$ as well as the neutral element of $\Z^d$.
We aim to give enough context so that this conflict causes no confusion.
\begin{thm}\label{thm: realize entropy smaller than lambda(k)}
    Let $G=\Z^d$, $k\geq 5$ and $0<h<\lambda(k)$.
    Then there exists a Toeplitz array $x\in (\Sigma_k)^G$ such that $(\ol{O_G}(x),G)$ is strictly ergodic and $\htop(\ol{O_G}(x),G) = h$.
    Moreover, the MEF of  $(\ol{O_G}(x),G)$ is given by the universal $G$-odometer.
\end{thm}
    Again, we shall divide the proof into several lemmas.
    We begin by defining two sequences $(p_n)_{n\in\N_0}$ and $(q_n)_{n\in\N_0}$ of integers that satisfy $p_{n+1} = p_n \cdot q_n$.
    Of course, for this it suffices to specify $p_0$ and all the $q_n$.
    Later we will use $p_n$ and $q_n$ to define diagonal matrices $P_n$ and $Q_n$.
    Let $M\geq 4$ be an integer such that
    \begin{equation}\label{eq: estimate xi-(d+1)log xi - 1}
        \xi - (d+1) \log(\xi) - 1 \geq 0 \text{ for all real numbers } \xi \geq M.
    \end{equation}
    Furthermore, let $N\in \N$ such that
    \begin{equation}\label{eq: choice of N}
    p_{N-1}^\prime h \geq 1, \quad q_{N-1}^\prime\geq M \quad \text{ and }\quad \frac{d\log(M)+3}{p_N^\prime} < \lambda(k)-h.
    \end{equation}
    Here, $(p_n^\prime)_{n\in\N_0}$, $(q_n^\prime)_{n\in\N_0}$, $(\lambda_n)_{n\in\N}$ are the sequences defined in Section \ref{subsection: three sequences}.
    We define $p_0 = p_0^\prime$ and $q_n = q_n^\prime$ for $0\leq n<N$.
    Note that this yields in particular $p_N = p_N^\prime$.
    It remains to define $q_{N+l}$ for $l\in \N_0$.
    \begin{lem}\label{lem: properties q_n}
        For every $l\in \N_0$ there exists an integer $q_{N+l}$ such that
        \begin{enumerate}[label=(\roman*-l)]
        \item $q_{N+l-1} <q_{N+l} \leq (q_{N+l-1} - 1)!$,
        \item $p_{N+l} \cdot h + 3 \leq \log q_{N+l} \leq p_{N+l}\cdot h+ d\log(M+l) + 3$,
        \item $(M+l)^d$ divides $q_{N+l}$.
        \item $M+l+1\leq q_{N+l}$.
    \end{enumerate}
    \end{lem}
    \begin{proof}
        The proof will be by induction on $l$.
        Note that (ii-0) is equivalent to
        $$e^3e^{p_{N} h} \leq q_{N} \leq M^d \cdot e^3e^{p_{N} h}.$$
        Since the interval $[e^3e^{p_{N} h}, M^d \cdot e^3e^{p_{N} h}]$ contains at least $M^d$ many consecutive integers, one can clearly find $q_N$ that satisfies (ii-0) and (iii-0).
        Moreover, (iv-0) will follow immediately from (i-0), since by \eqref{eq: choice of N} we have $q_N > q_{N-1} = q_{N-1}^\prime  \geq M$,
        so that the fact that $q_N$ is an integer yields $q_N \geq M+1$.
        If we prove that (ii-0) implies (i-0), then the proof of the induction start with $l=0$ will be complete.
        This implication will follow from the two inequalities
        $$e^{p_{N}h + 3} > q_{N-1}\quad\text{ and }\quad\log((q_{N-1} -1)!) \geq p_{N}\cdot h + d\log(M) + 3.$$
        Observe for the first inequality, that by \eqref{eq: choice of N} we have $p_{N-1}^\prime \cdot h \geq 1$.
        Hence,
        $$ e^{p_{N}\cdot h + 3} > e^{q_{N-1}(p_{N-1}\cdot h)} = e^{q_{N-1}(p_{N-1}^\prime \cdot h)} \geq e^{q_{N-1}} \geq  q_{N-1}.$$
        For the second inequality, observe that
        \begin{align*}
        \log((q_{N -1}-1)!) &= \log q_{N}^\prime
        = p_{N}^\prime \lambda_{N}
        \geq p_{N}^\prime \cdot\lambda(k)
        \overset{\eqref{eq: choice of N}}{\geq} p_{N}^\prime\cdot  h+d\log(M) + 3.
        \end{align*}
        Due to $p_N^\prime = p_N$, this concludes the proof of the case $l=0$.\\
        Now suppose that we have found an integer $q_{N+l}$ which satisfies properties (i-l) - (iv-l).
        Analogously to the base case, we see that we can find $q_{N+l+1}$ that satisfies (ii-(l+1)) and (iii-(l+1)).
        Moreover we see that (i-(l+1)) implies (iv-(l+1)), as
        $$ q_{N+l+1} \overset{\text{(i-(l+1))}}> q_{N+l} \overset{\text{(iv-l)}}{\geq} M + l + 1.$$
        Similarly to above, the induction step will be finished, if we show that
        $$e^{p_{N+l+1} h + 3}> q_{N+l}\quad \text{ and }\quad\log((q_{N+l}-1)!) \geq p_{N+l+1}\cdot h + d\log(M+l+1)+3.$$
        For the first inequality, recall that the $p_n$ are monotonically increasing, so that we have $p_{N+l} \cdot h \geq p_{N-1} \cdot h = p_{N-1}^\prime \cdot h \geq 1$.
        Thus,
        $$ e^{p_{N+l+1}\cdot h + 3} > e^{q_{N+l}(p_{N+l}\cdot h)} \geq e^{q_{N+1}} \geq  q_{N+l}.$$
        For the second inequality, we observe that
    \begin{align*}
        \log((q_{N+l}-1)!) &\geq (q_{N+l} - 1) (\log(q_{N+l} - 1) -1)     &&\eqref{eq: estimate factorial}\\
        &\geq (q_{N+l}-1) (\log(q_{N+l})-2)   &&(\text{mean val.\ thm.\ on log})\\
        &= q_{N+l} \log(q_{N+l}) - 2q_{N+l} - \log(q_{N+l}) + 2\\
        &\geq  q_{N+l} \cdot p_{N+l}\cdot h\cdot + q_{N+l} - \log(q_{N+l}) + 2  &&\text{(ii-l)}\\
        &= p_{N+l+1} \cdot h+ q_{N+l} - \log(q_{N+l}) + 2.
    \end{align*}
    Now note that by (iv-l) we have $q_{N+l} \geq M+l+1$ and in particular $q_{N+l} \geq M$.
    Thus,
    \begin{align*}
        q_{N+l} - \log(q_{N+l}) + 2 - (d\log(M+l+1) + 3) \geq q_{N+l} - (d+1)\log(q_{N+l}) - 1  \geq 0,
    \end{align*}
    where the last inequality follows by \eqref{eq: estimate xi-(d+1)log xi - 1}.
    This proves the desired inequality and hence the lemma.
    \end{proof}
    Now we can define diagonal matrices $P_n$ and $Q_n$ such that $P_{n+1} = P_n \cdot Q_n$.
    We define $P_1 = \diag(p_1,1,\ldots,1)$ and for $n<N$ we set $Q_n = \diag(q_n,1,\ldots,1)$.
    Moreover, we define for $l\in \N_0$
    $$ Q_{N+l} = \diag\kl \frac{q_{N+l}}{(M+l)^{d-1}},M+l,M+l,\ldots,M+l\kr.$$
    By Lem.\ \ref{lem: properties q_n} (iii-l), we see that $Q_{N+l} \in \Z^{d\times d}$ and moreover $M+l$ divides every diagonal entry of $Q_{N+l}$ and thus also of $P_{N+l+1}$.\\
    This uniquely determines the scale $(P_n)_{n\in\N}$ with increment $(Q_n)_{n\in\N}$. We denote by $\Gamma_n = P_n \Z^d$ the subgroups associated to $P_n$.
    One sees quickly that all conditions of Lem.\ \ref{lem: shifting cuboids for better fundamental domain} are fulfilled, so that this lemma provides us with a suitable sequence $(D_n)_{n\in\N}$ of fundamental domains of $G/\Gamma_n$.
    It is easy to see that $\# D_n = p_n$ and $\# (D_{n+1}\cap \Gamma_n) = q_n$.
    Choose for $n\in \N$ a subset $\mathcal{S}_n$ of
    $$ \setv{\sigma: (D_{n}\cap \Gamma_{n-1})\setminus\set{0}\to \set{2,\ldots,q_{n-1}}}{\sigma \text{ is bijective}}$$
    of cardinality $\# \mathcal{S}_n = q_n$ such that
    \begin{equation}\label{eq: every block is send to every position}
       \forall \gamma\in (D_n \cap \Gamma_{n-1})\setminus\set{0}: \forall j \in \set{2,\ldots,q_{n-1}}: \exists \sigma\in \mathcal{S}_n:\: \sigma(\gamma)= j.
    \end{equation}
    Note that this can be done for $n\geq N$ due to (i-l) (with $l:= n-N\geq 0$).
    For $n<N$ it is also clear since there is no actual choice for $\mathcal{S}_n$ because the condition
    $$ \# \mathcal{S}_n = q_n = q_n^\prime = (q_{n-1}^\prime-1)! = (q_{n-1}-1)!$$
    forces to take the whole set
    $$ \mathcal{S}_n = \setv{\sigma: (D_{n}\cap \Gamma_{n-1})\setminus\set{0}\to \set{2,\ldots,q_{n-1}}}{\sigma \text{ is bijective}}.$$
    This choice clearly satisfies \eqref{eq: every block is send to every position}.
    Next, we define sets $\mathcal{C}_n$ of $D_n$-blocks with cardinality $\# \mathcal{C}_n = q_n$.
    Recall that we let $\Gamma_0 = G$ and $D_0 = \set{0} \tm G$.
    For $j\in \set{1,\ldots,k}$, let
    $$ C_0^{(j)}: \set{0}=D_0 \to \Sigma_k,\: 0 \mapsto j-1,$$
    so that $\mathcal{C}_0 := \setv{C_0^{(j)}}{j\in\set{1,\ldots,k}}$ has cardinality $k=q_0^\prime = q_0$.
    Now let $n\in\N_0$ and assume that we have already defined the set $\mathcal{C}_n = \setv{C_n^{(j)}}{j\in \set{1,\ldots,q_n}}$.
    We can define for every $\sigma\in \mathcal{S}_{n+1}$ a $D_{n+1}$-block $C_{n+1}^\sigma$ by setting
    \begin{equation}\label{eq: def C_{n+1}}
    C_{n+1}^\sigma (d_n + \gamma_n) = \begin{cases}
        C_n^{(1)}(d_n)                      &, \gamma_n =0\\
        C_n^{(\sigma(\gamma_n))}(d_n)       &, \gamma_n\neq 0
    \end{cases},
    \end{equation}
    where $d_n \in D_n$ and $\gamma_n\in D_{n+1}\cap \Gamma_n$.
    This is well-defined due to Lem.\ \ref{lem: shifting cuboids for better fundamental domain} (iv).
    Now we can let $\mathcal{C}_{n+1} = \setv{C_{n+1}^\sigma}{\sigma\in\mathcal{S}_{n+1}}$ and give it some arbitrary ordering $\mathcal{C}_{n+1} = \setv{C_{n+1}^{(j)}}{j\in \set{1,\ldots,q_{n+1}}}$.\\
    Let the maps $\theta_i: G \to D_i\cap \Gamma_{i-1}$ ($i\in \N$) be defined as in Def.\ \ref{def: theta maps}.
    \begin{lem}\label{lem: going backwards for the C_n}
        Let $m>n\geq 0$.
        For every $C\in \mathcal{C}_m$ and $g\in D_m$ there exists $j\in \set{1,\ldots,q_n}$ such that
        $$ C(g) = C_n^{(j)}\kl \sum_{i=1}^n \theta_i(g)\kr.$$
        If $\theta_{n+1}(g) = 0$, then $j=1$.
    \end{lem}
    The proof is a straightforward induction on $m$.
    Now we can define the Toeplitz array $x$.
    We will do this by identifying $x$ as the limit of a sequence of $\Gamma_n$-periodic arrays $y_n$ over the alphabet $\Sigma_k\cup \set{*}$.
    Let $g\in G$ and define $y_0(g) = *$. For $n\geq 1$, let
    $$ y_n(g) = \begin{cases}
        C_{n-1}^{(1)}\kl \sum_{i=1}^{n-1}\theta_i(g)\kr            &, \theta_n(g) = 0\\
        y_{n-1}\kl \sum_{i=1}^{n-1}\theta_i(g)\kr           &, \theta_n(g)\neq 0
    \end{cases}.$$
    By Lem.\ \ref{lem: going backwards for the C_n}, it follows that $y_m\vert_{D_n} = C_n^{(1)}$ for all $m>n$.
    As $G= \bigcup_{n\in \N} D_n$, we see that $(y_n)_{n\in\N_0}$ converges to some array $x\in (\Sigma_k)^G$.

    \begin{rem}\label{rem: explain construction}
        One can understand this construction as follows:
        We consider the symbols $*$ as "holes".
        In each step of the construction we take some holes, "fill" them with symbols from $\Sigma_k$ and repeat this filling with some period.
        As soon as a hole has been filled in, it will never be changed in the future again.
        We also make sure that every hole is filled eventually.
        This method has been used for numerous constructions of Toeplitz $\Z$-sequences.
        See for example \cite{Downarowicz2005Survey}, \cite{IwanikLacroix1994NonRegular} and \cite{Williams1984ToeplitzFlows}.
    \end{rem}

    \begin{lem}
        Let $n\in\N$ and $g\in G$.
        \begin{enumerate}[label=(\roman*)]
            \item We have $y_n(g)=*$ iff $\theta_i(g)\neq 0$ for all $i\in\set{1,\ldots,n}$.
            \item Let $y_n(g)\neq *$ and let $n\geq 1$ be chosen minimal with this property.
            Then, we have $y_m(g)=y_n(g)$ for all $m>n$.
        \end{enumerate}
    \end{lem}
    \begin{proof}
        (i) This follows immediately from the definition of $y_n$.\\
        (ii) If $n$ is minimal with $y_n(g)\neq *$, then by (i) it follows that $\theta_n(g) = 0$ and $\theta_i(g)\neq 0$ for all $i<n$.
        Furthermore, we obtain $y_n(g) = C_{n-1}^{(1)}\kl \sum_{i=1}^{n-1}\theta_i(g)\kr.$
        Now let $m>n$.
        If $\theta_m(g)\neq 0$, then $y_m(g) = y_{m-1}(g)$ and we can argue by induction on $m$.
        On the other hand, if $\theta_m(g)=0$, then
        \[ y_m(g) = C_{m-1}^{(1)} \kl \sum_{i=1}^{m-1} \theta_i(g)\kr \overset{\text{Lem.\ \ref{lem: going backwards for the C_n}}}{=} C_{n-1}^{(1)}\kl \sum_{i=1}^{n-1}\theta_i(g) \kr = y_n(g). \qedhere\]
    \end{proof}
    It follows from the previous lemma that if $y_n(g)\neq *$, then $g\in \Per(x,\Gamma_n)$, so that
    \begin{equation}\label{eq: periodic positions of x inclusion}
        \Per(x,\Gamma_n) \mt \setv{g\in G}{\theta_i(g) = 0 \text{ for some } i\in\set{1,\ldots,n}}.
    \end{equation}
    In particular, $\Per(x,\Gamma_n)\mt D_{n-1}$, hence $x$ is a Toeplitz array.
    We also conclude that $0 \in \Per(x,\Gamma_n)$ for every $n\in\N$.
    Because of $y_1(0) = C_0^{(1)}(0) = 0$, we obtain
    \begin{equation}\label{eq: 0 is Gamma_n-periodic}
        0 \in \Per(x,\Gamma_n,0) \text{ for every } n\in \N.
    \end{equation}
    We now aim to prove the converse inclusion of \eqref{eq: periodic positions of x inclusion}.
    It turns out that we are even able to prove a stronger statement:
    \begin{lem}\label{lem: every letter can be found on coset}
        Let $g\in G$ with $\theta_i(g)\neq 0$ for all $i\in \set{1,\ldots,n}$ and let $\alpha\in \Sigma_k\setminus \set{0} = \set{1,\ldots,k-1}$.
        Then, there exists $\gamma\in \Gamma_n$ such that $x(g+\gamma) = \alpha$.
    \end{lem}
    \begin{proof}
    We can assume w.l.o.g.\ that $g\in D_n$.
    Let also $\alpha \in \Sigma_k\setminus\set{0} = \set{1,\ldots,k-1}$.
    We claim that there exists $j\in\set{2,\ldots,q_n}$ such that $C_n^{(j)}(g)=\alpha$.
    For $n=0$ this is fulfilled by definition of $C_0^{(2)},\ldots,C_0^{(k)}$.
    Now let $n\geq 1$, $g\in D_{n}$ and let $\theta_i(g) \neq 0$ for all $i\in\set{1,\ldots,n}$.
    We can write $g = h+ \theta_{n}(g)$ with $h\in D_{n-1}$.
    Since $\theta_i(h)=\theta_i(g)\neq 0$ for all $i\in \set{1,\ldots,n-1}$, we obtain by induction some $l\in \set{2,\ldots,q_{n-1}}$ such that $C_{n-1}^{(l)}(h) = \alpha$.
    By \eqref{eq: every block is send to every position} and the fact that $\theta_{n}(g)\neq 0$, there exists $j\in\set{2,\ldots,q_{n}}$ such that $C_{n}^{(j)}(g) = C_{n-1}^{(l)}(h) = \alpha$, proving the claim.\\
    Furthermore, by definition of $C_{n+1}^{(1)}$, there exists $\gamma\in D_{n+1}\cap \Gamma_n$ such that $C_{n+1}^{(1)}(g+\gamma) = C_n^{(j)}(g)$.
    Together with the claim, we obtain
    \[
        x(g+\gamma) = y_{n+2}(g+\gamma) =  C_{n+1}^{(1)}(g+\gamma) = C_n^{(j)}(g) = \alpha. \qedhere
    \]
    \end{proof}
    In particular, since $\alpha$ can be chosen arbitrarily and $\Sigma_k\setminus\set{0}$ contains at least two elements (because of $k\geq 5$), it follows that $g\notin\Per(x,\Gamma_n)$ for any $g$ with $\theta_i(g)\neq 0$ for all $i\in \set{1,\ldots,n}$.
    Thus, we have shown
    \begin{equation}\label{eq: periodic positions of x}
        \Per(x,\Gamma_n) = \setv{g\in G}{\theta_i(g) = 0 \text{ for some } i \in \set{1,\ldots,n}}.
    \end{equation}
    Therefore, Lem.\ \ref{lem: property * implies essential period} and Lem.\ \ref{lem: properties theta maps} yield that every $\Gamma_n$ is an essential period of $x$.
    The fact that the universal $G$-odometer is the MEF of $(\ol{O_G}(x),G)$ follows from Prop.\ \ref{prop: MEF of Toeplitz} and Cor.\ \ref{cor: universal odo. on Z^d} as well as the construction of the matrices $P_n$ ($n\in\N$).\\
    Clearly we have $W_{D_n}(x)= \mathcal{C}_n$, so that $\# W_{D_n}(x)=q_n$ for all $n\in\N$.
    Moreover, property (ii-l) yields for all $l\in \N_0$ that
    $$ h+ \frac{3}{p_{N+l}} \leq \frac{\log q_{N+l}}{p_{N+l}} \leq h + \frac{d \log(M+l) + 3}{p_{N+l}}.$$
    Since clearly $p_{N+l}\geq l$, it follows with Prop.\ \ref{prop: entropy Toeplitz} that
    $$ \htop(\ol{O_G}(x),G) = \lim_{n\to \infty} \frac{\log \# W_{D_n}(x)}{\# D_n} = \lim_{n\to\infty} \frac{\log q_n}{p_n} = h.$$
    Lastly, \eqref{eq: def C_{n+1}} implies $\mathrm{ap}(B,C) = 1/q_n$ for all $B\in W_{D_n}(x)$ and $C\in W_{D_{n+1}}(x)$, so that strict ergodicity follows by Cor.\ \ref{cor: sufficient crit. unique ergod. Toeplitz}.
    This finishes the proof of Thm.\ \ref{thm: realize entropy smaller than lambda(k)}. \qed\\

Now we are in position to prove the main result.

\begin{proof}[Proof of Theorem A]
    Let $k\geq 2$ and $0<h<\log k$ be arbitrary.
    We choose $s\in \N$ such that $\sqrt[d]{s}\in \Z$, $k^{s}\geq 5$ and $s(\log k - h ) \geq 5$.
    By Lem.\ \ref{lem: properties p_n, lambda_n}, we deduce that
    $$ \lambda(k^{s}) > \log k^{s} - 5 = s \log k - 5 \geq s h.$$
    By Thm.\ \ref{thm: realize entropy smaller than lambda(k)}, we find a Toeplitz array $x \in (\Sigma_{k^{s}})^G$ such that $(\ol{O_G}(x),G)$ is strictly ergodic and $\htop(\ol{O_G}(x),G)= s h$.
    Moreover, the proof of this theorem yields that a period structure $(\Gamma_n)_{n\in\N}$ of $x$ is given by subgroups $\Gamma_n = P_n \Z^d$ with diagonal matrices $P_n$ with positive diagonal entries.
    The scale $(P_n)_{n\in\N}$ satisfies the conditions of Cor.\ \ref{cor: universal odo. on Z^d}.
    Let $P = \diag(\sqrt[d]{s},\ldots,\sqrt[d]{s})\in \Z^{d\times d}$ and
    $$F = P([0,1)^d)\cap \Z^d = \setv{g=(g_1,\ldots,g_d)\in \Z^d}{\forall i\in \set{1,\ldots,d}: \: 0\leq g_i < \sqrt[d]{s}}.$$
    Note that $\# F = s$.
    Take some bijective map $S: \Sigma_{k^{s}} \to (\Sigma_k)^F$.
    (We will later impose some additional conditions on $S$.)
    Now, Prop.\ \ref{prop: entropy and strict ergod. under substitution} immediately yields that for $y:= x^S\in (\Sigma_k)^G$ we have $\htop(\ol{O_G}(y),G) = \htop(\ol{O_G}(x),G)/\# F = (sh)/s = h$ and the Toeplitz subshift $(\ol{O_G}(y),G)$ is strictly ergodic.
    Obviously, the matrices $P\cdot P_n$ still satisfy the conditions of Cor.\ \ref{cor: universal odo. on Z^d}, so that the associated odometer is also the universal $\Z^d$-odometer.
    Hence, the proof will be finished if we show that the sequence $(\Delta_n)_{n\in\N}$ with  $\Delta_n := (P\cdot P_n)\Z^d = P\Gamma_n$ is a period structure of $y$.\\
    In order to do so, we aim to use Lem.\ \ref{lem: property * implies essential period}. Therefore, we seek to determine the set $\Per(y,\Delta_n)$.
    First recall from the discussion prior to Rem.\ \ref{rem: period structure after substitution} that $F_n = F + PD_n$ is a fundamental domain of $G/\Delta_n$ which satisfies all properties (i) - (v) of Lem.\ \ref{lem: shifting cuboids for better fundamental domain}.
    Let $\theta_i^{\Gamma}: G\to D_i \cap \Gamma_{i-1}$ ($i\in \N$) denote the maps given by Def.\ \ref{def: theta maps} applied to $(\Gamma_n)_{n\in\N}$ and $(D_n)_{n\in\N}$.
    We wish to apply the same definition to obtain maps $\theta_i^{\Delta}$.
    However, observe that $\Delta_0 = P\Gamma_0 = P\Z^d \neq \Z^d$ whereas Def.\ \ref{def: theta maps} requires the convention $\Delta_0 = G = \Z^d$.
    This forces us to add additional members to the sequences, more specifically defining $F_0 = F$ and $\Delta_{-1} = G$.
    With these precautions, applying Def.\ \ref{def: theta maps} to $(\Delta_n)_{n\in\N_0}$ and $(F_n)_{n\in\N_0}$ with the index starting at $i=0$ yields maps $\theta_i^\Delta: G \to F_i \cap \Delta_{i-1}$ ($i\in \N_0$).
    Let $g\in F_n$ and write $g = f + Pd$ with $f\in F$ and $d\in D_n$.
    It is easy to see that
    $$ \theta_i^{\Delta}(g) = \begin{cases}
        f                   &, i = 0\\
        P\theta_i^{\Gamma}(d) &, i>0
    \end{cases}. $$
    Assume that $\theta_i^{\Delta}(g)\neq 0$ for all $i\in\set{1,\ldots,n}$.
    Hence, $\theta_i^{\Gamma}(d) \neq 0$ for all $i\in \set{1,\ldots,n}$.
    By Lem.\ \ref{lem: every letter can be found on coset}, for every $\alpha\in \Sigma_{k^{s}}\setminus\set{0}$ there exists $\gamma_\alpha\in \Gamma_n$ such that $x(d+\gamma_\alpha) = \alpha$.
    Let $\beta = x(d)\in \Sigma_{k^{s}}$.
    As we can assume $s\geq 2$, it follows that $k^{s}-1 > k^{s-1}$. So by bijectivity of $S: \Sigma_{k^{s}}\to (\Sigma_k)^F$, there exists $\alpha\in \Sigma_{k^{s}} \setminus\set{0}$ such that $[S(\alpha)](f)\neq [S(\beta)](f)$.
    Now we can compute that
    \begin{align*}
        y(g+P\gamma_\alpha)= x^S(f+P(d+\gamma_\alpha)) =[S(\alpha)](f)\neq [S(\beta)](f) = x^S(f+Pd) =y(g).
    \end{align*}
    Since $P\gamma_\alpha\in \Delta_n$, it follows that $g\notin \Per(y,\Delta_n)$.
    Conversely, if $\theta_i^{\Delta}(g)=0$ for some $i\in\set{1,\ldots,n}$, then by \eqref{eq: periodic positions of x} it easily follows that $g\in \Per(y,\Delta_n)$.
    Hence,
    \begin{equation}\label{eq: periodic positions of y}
        \Per(y,\Delta_n)=\setv{g\in G}{\theta_i^{\Delta}(g) = 0 \text{ for some } i\in \set{1,\ldots, n}}.
    \end{equation}
    This implies in particular that $G = \bigcup_{n\in\N}\Per(y,\Delta_n)$ and $\Per(y,\Delta_n)\neq \emptyset$ for every $n\in\N$.
    Thus, it remains to show for arbitrary $n\in\N$ and $h\in G$ that the implication
    $$(\forall \alpha\in\Sigma_k: \:\Per(y,\Delta_n,\alpha)\tm\Per(y,\Delta_n,\alpha)-h) \implies h\in \Delta_n$$
    holds, i.e.\ that $\Delta_n$ is essential.
    Let $h\in G$ such that $\Per(y,\Delta_n,\alpha)\tm \Per(y,\Delta_n,\alpha)-h$ for all $\alpha \in \Sigma_k$ and assume for a contradiction that $h\notin \Delta_n$.
    We immediately deduce using \eqref{eq: periodic positions of y} and Lem.\ \ref{lem: property * implies essential period} that $\theta_i^{\Delta}(h) = 0$ for all $i\in \set{1,\ldots,n}$.
    Hence, we can assume w.l.o.g.\ that $h = \theta_0^{\Delta}(h) \in F\setminus \set{0}$.
    Note that we can make additional assumptions on the substitution $S$ in order to make the proof more convenient and direct.
    Therefore, we specify the block $S(0)\in(\Sigma_k)^F$ by assuming for $f\in F$ that
    $$ [S(0)](f) = \begin{cases}
        0   &, f=0\\
        1   &, f\neq 0
    \end{cases}.$$
    (However, it should be noted that it is possible to arrive at the desired contradiction without any assumptions on $S$ other than bijectivity.)
    Recall that $0\in \Per(x,\Gamma_n,0)$ by \eqref{eq: 0 is Gamma_n-periodic}, so that $0\in \Per(y,\Delta_n)$.
    Moreover, due to the choice of $S$ we have $y(0) = [S(0)](0) = 0$, so that $0\in \Per(y,\Delta_n,0)$.
    However, it also follows that
    $$y(0+h) =[S(x(0))](h) = [S(0)](h) = 1,$$
    hence certainly $0+h \notin \Per(x,\Delta_n,0)$, which is our desired contradiction.
    We have thus shown that $(\Delta_n)_{n\in\N}$ is indeed a period structure of $y$, hence the MEF of $(\ol{O_G}(y),G)$ is the universal $\Z^d$-odometer.
    This concludes the proof of Theorem A.
\end{proof}

\subsection{Consequences}\label{subsection: consequences}
We want to point out an immediate consequence of Theorem A.
Namely, combining it with a result due to Cortez and Gomez allows us to additionally realize any (non-periodic) odometer as the MEF of a strictly ergodic almost 1-1 extension with any prescribed entropy.
\begin{thm}[{\cite[Thm.\ 1]{CortezGomez2024Amenability}}]\label{thm: find almost 1-1 between odometer and extension}
    Let $(\overleftarrow{G},\Z^d)$ be a non-periodic $\Z^d$-odometer and let $(X,\Z^d)$ be a minimal extension of $(\overleftarrow{G},\Z^d)$ with factor map $\pi_X: X \to \overleftarrow{G}$.
    Then, there exists a minimal topological dynamical system $(Y,\Z^d)$ and an almost 1-1 factor map $\pi_Y: Y\to \overleftarrow{G}$ such that
    \begin{enumerate}[label=(\roman*)]
        \item There exists a factor map $\varphi: X\to Y$ such that $\varphi$ is 1-1 on $\pi_X^{-1}(Z_1)$ where $Z_1\tm \overleftarrow{G}$ is a $\Z^d$-invariant set with full Haar measure.
        \item $\pi_X = \pi_Y \circ \varphi$.
        \item The map $\varphi$ induces an affine bijection between the spaces of invariant Borel probability measures on $(X,\Z^d)$ and $(Y,\Z^d)$.
    \end{enumerate}
    If $X$ is a Cantor set, then so is $Y$.
\end{thm}
Note that due to the Variational Principle, the map $\varphi$ preserves topological entropy.
\begin{cor}\label{cor: find almost 1-1 between odometer and extension (strictly ergodic)}
    Let $(\overleftarrow{G},\Z^d)$ be a non-periodic $\Z^d$-odometer and let $(X,\Z^d)$ be a strictly ergodic extension of $(\overleftarrow{G},\Z^d)$.
    Then, there exists a strictly ergodic almost 1-1 extension $(Y,\Z^d)$ of $(\overleftarrow{G},\Z^d)$ such that $\htop(Y,\Z^d) = \htop(X,\Z^d)$.
    If $X$ is a Cantor set, then so is $Y$
\end{cor}

Hence, together with Theorem A and the fact that any odometer is a factor of the universal odometer we obtain the following:

\begin{namedtheorem}[Theorem A']
    Let $h>0$ and $(\overleftarrow{G},\Z^d)$ be a non-periodic $\Z^d$-odometer.
    There exists a strictly ergodic topological dynamical system $(Y,\Z^d)$ on a Cantor set $Y$, which is an almost 1-1 extension of $(\overleftarrow{G},\Z^d)$ and which satisfies $\htop(Y,\Z^d)= h$.
\end{namedtheorem}

\bibliographystyle{alpha}

\end{document}